\documentclass[review]{elsarticle}

\usepackage[dvipdfmx]{hyperref}
\journal{Elsevier}

\usepackage{amsmath}
\usepackage{amssymb}
\usepackage{algorithm}
\usepackage{algorithmic}
\newcommand{\bs}[1]{\boldsymbol{#1}}

\newcommand{\figref}[1]{{Fig. \ref{#1}}}
\newcommand{\equref}[1]{{Eq. (\ref{#1})}}

\begin{document}

\begin{frontmatter}

\title{A topology optimisation of acoustic devices based on the frequency response estimation with the Pad\'{e} approximation}
% \tnotetext[mytitlenote]{Fully documented templates are available in the elsarticle package on \href{http://www.ctan.org/tex-archive/macros/latex/contrib/elsarticle}{CTAN}.}

% \address{Radarweg 29, Amsterdam}
% \fntext[myfootnote]{Since 1880.}

\author[nagoya]{Yuta Honshuku}

\author[keio]{Hiroshi Isakari\corref{mycorrespondingauthor}}
\cortext[mycorrespondingauthor]{Corresponding author}
\ead{isakari@sd.keio.ac.jp}

%\author[nagoya]{Toru Takahashi}
%\author[nagoya]{Toshiro Matsumoto}

\address[nagoya]{Nagoya University, Furo-cho, Chikusa-ku, Nagoya, Aichi, Japan}
\address[keio]{Keio University, 3-14-1, Hiyoshi, Kohoku-ku, Yokohama, Kanagawa, Japan}

\begin{abstract}
We propose a topology optimisation of acoustic devices that work in a certain bandwidth. To achieve this, we define the objective function as the frequency-averaged sound intensity at given observation points, which is represented by a frequency integral over a given frequency band. It is, however, prohibitively expensive to evaluate such an integral naively by a quadrature. We thus estimate the frequency response by the Pad\'{e} approximation and integrate the approximated function to obtain the objective function. The corresponding topological derivative is derived with the help of the adjoint variable method and chain rule. It is shown that the objective and its sensitivity can be evaluated semi-analytically. We present efficient numerical procedures to compute them and incorporate them into a topology optimisation based on the level-set method. We confirm the validity and effectiveness of the present method through some numerical examples.
\end{abstract}

\begin{keyword}
\texttt{topology optimisation, fast frequency sweep, Pad\'{e} approximation, acoustic device, boundary element method, topological derivative}
\MSC[2010] 00-01\sep  99-00
\end{keyword}

\end{frontmatter}

%\linenumbers

\section{Introduction} \label{intro_section}
Controlling sound waves in the desired manner is one of the most significant tasks in engineering. Acoustic devices to manipulate sound waves are currently designed with some trial-and-error, while the recent rapid development in computer-aided design has enabled us to exploit structural optimisation. The structural optimisation methods are classified into sizing~\cite{prager1974note}, shape~\cite{takahashi2021shape}, and topology optimisations~\cite{bendsoe1988generating, amstutz2006new, isakari2017Bspline}, the last of which has the highest degrees of design freedom. The topology optimisation regards structural optimisation as material distribution optimisation, which allows topological changes (e.g. creation of new holes and/or inclusions) in its process. Various researches on topology optimisation for acoustic devices are found in the literature, e.g.~\cite{wadbro2006topology, kook2013acoustical}. Typically, the existing methods optimise the acoustic response to an incident wave of a single frequency. The frequency response, however, highly depends on the excitation frequency. A device optimised to a target frequency is specialised to it, and thus its performance may fall due to even a slight frequency fluctuation. When the frequency variation is expected in practical use, such a structure is not desirable. The acoustic device should exhibit excellent performance in the whole expected frequency range.

One possible approach to realise such a wide-band device is robust topology optimisation~\cite{sato2020robust, qin2021robust}. In the robust topology optimisation, the (angular) frequency is regarded as a stochastic variable following the normal distribution. The frequency response becomes a stochastic variable accordingly. With these settings, one can obtain an acoustic device robust to the frequency fluctuation by letting the objective function be, for example, the weighted sum of the original objective and its standard deviation. Since it is time-consuming to evaluate the statistic quantity by the Monte-Carlo method, the deviation is often approximated through the Taylor expansion of the original objective at the mean angular frequency. We have pointed out that the high-order approximation is sometimes necessary to achieve a wide-band topology optimisation in acoustics and developed a low-cost method to evaluate the angular frequency derivatives up to an arbitrary order which is required in such a high-order approximation~\cite{qin2021robust}.

Another approach to realise a wide-band optimisation may define the objective function $J$ as the frequency average of the original objective $f(\omega)$ over a given frequency range $[\omega_1,\omega_2]$ (target band), which is formulated as follows:
\begin{equation}
 J=\frac{1}{\omega_2-\omega_1}\int_{\omega_1}^{\omega_2} f(\omega) \mathrm{d}\omega.
  \label{frequency_average}
\end{equation}
Since $f(\omega)$ is generally unknown, we need some ingenuity to evaluate the integral in \equref{frequency_average}. The simplest solution may utilise a numerical quadrature. This is, however, computationally demanding because we need to repeat the response analysis $N$ times for various excitation frequency $\omega$ to evaluate the integrand if we use a quadrature with $N$ integral points. Another solution is to replace the integrand $f(\omega)$ with its approximant. Jensen~\cite{jensen2007topology} proposed such a method based on the Pad\'{e} approximation. In \cite{jensen2007topology}, $f(\omega)$ is defined as:
\begin{equation}
 f(\omega)=\sum_{i=1}^{N_\mathrm{obs}}L_i|u(\omega; \boldsymbol{x}_{\mathrm{obs}}^i)|^2,
\end{equation}
where $\bs{x}_\mathrm{obs}^i$, $L_i$, and $u$ represent the observation points, given weighting coefficient, and complex amplitude, respectively. In the evaluation of $J$ in \equref{frequency_average}, $u(\omega; \bs{x}^i_\mathrm{obs})$ is first approximated by an extended-version of the Pad\'{e} approximation which can approximate $u(\omega;\bs{x}_\mathrm{obs}^i)$ at all the observation points $\bs{x}_\mathrm{obs}^i~(i=1,\cdots, N_\mathrm{obs})$ simultaneously. The approximant for $u(\omega; \bs{x}_\mathrm{obs}^i)$ is then used to form an approximation of $f(\omega)$, which is integrated by a numerical quadrature to evaluate (the approximation of) $J$. Note that the frequency range where such an approximation for $f$ holds is limited to a certain range around the approximation centre $\omega_0$. In the case that a single Pad\'{e} approximation cannot cover the entire target band $[\omega_1, \omega_2]$, the band is divided into several subintervals in each of which the approximation is applied. In \cite{jensen2007topology}, the number of the subintervals is given a priori.

The Pad\'{e} approximation uses the following rational function:
\begin{equation}
 g_{[M,N]}(z):=\frac{\displaystyle{\sum_{i=0}^M}p_i(z-z_0)^i}{\displaystyle{\sum_{i=0}^N}q_i(z-z_0)^i}
\end{equation}
to approximate a function $g(z)$ around $z_0$, in which $g: \mathbb{C}\rightarrow\mathbb{C}$ is generally a complex-valued function, and $p_i~(i=0,\cdots,M)$ and $q_i~(i=0,\cdots,N)\in\mathbb{C}$ are the coefficients. We call the rational function $g_{[M,N]}$ as ``$[M,N]$ Pad\'{e} approximation'' of $g$ in this paper. The coefficients $p_i$ and $q_i$ are determined such that the approximated rational polynomial $g_{[M,N]}(z)$ has the same Taylor expansion at $z_0$ as the original function $g(z)$ up to $(M+N)$-th order. The most important advantage of the Pad\'{e} approximation is that it can beyond the Taylor expansion's radius of convergence. The frequency response function $u(\omega;\bs{x}_\mathrm{obs}^i)$ may have some poles corresponding to the eigenfrequencies, and the radius of convergence is determined by the distance from the expansion centre to the nearest pole. The frequency range where the Taylor expansion is valid can thus be very narrow. This is exactly why the rational polynomial instead of the standard polynomial is used in \cite{jensen2007topology}. It is shown that the strategy may work even for a practical design in  automobile industry~\cite{kook2013acoustical}.

This paper aims to extend the Pad\'{e}-based topology optimisation to manipulate acoustic scattering in an unbounded region, which is of great interest in engineering because many acoustic devices (e.g. sound-proofing wall) operate in the infinite domain. We adopt the boundary element method (BEM) instead of the standard finite element method (FEM) in the objective and sensitivity computations. This choice is motivated by the fact that BEM requires fewer elements than FEM in the scattering analysis. Besides, we can deal with the boundary condition at the infinity in a strict manner by BEM without any artificial boundary. In addition to the extension, we shall address improving the original Pad\'{e}-based topology optimisation in two ways. Firstly, we propose to avoid the numerical integration to evaluate the frequency integral of the Pad\'{e}-approximated response. To use the numerical quadrature for the rational polynomial, one needs to carefully choose the quadrature rule since the rational polynomial may have some poles. Otherwise, the result can be inaccurate. It is also difficult to set in advance the appropriate number of integral points. We propose a drastic solution to these issues, i.e. derive the analytical expression for the objective function written by the integral of the rational polynomial approximation. We also present the analytical expression for its topological derivative. Furthermore, we provide an efficient numerical algorithm to compute the ingredients of the quantities. Secondly, we propose a simple adaptive strategy to subdivide the integral range in the case that we cannot accurately evaluate the objective with a single Pad\'{e} approximation.

The rest of this paper is organised as follows. In Section \ref{statement_of_the_problem}, we formulate the boundary value problem governing the scattering problem and the related structural optimisation problem. In Section \ref{objective_func_section}, we derive the approximations of the frequency average based on the Pad\'{e} approximation. We propose two approximations and compare their performances numerically in Section \ref{sec:comparison_of_the_two_methods}. In Sections \ref{freq_der_section} and \ref{a_fast_method_to_compute_the_rhs_of} we review the way to evaluate the angular frequency derivatives, which is required to obtain the Pad\'{e} approximation. In Section \ref{sensitivity_section}, we derive the topological derivative of the objective function. In Section \ref{range_division}, we propose a simple method to divide the target band into the appropriate number of subintervals. In Section \ref{example_section}, we show some validation results and examples of optimised acoustic designs. Finally, we conclude this paper in Section \ref{conclusion_section}.

\section{Statement of the problem}
\label{statement_of_the_problem}
In this section, we define the boundary value problem and the structural optimisation problem of interest.

\begin{figure}
 \centering
 \scalebox{0.8}{\includegraphics{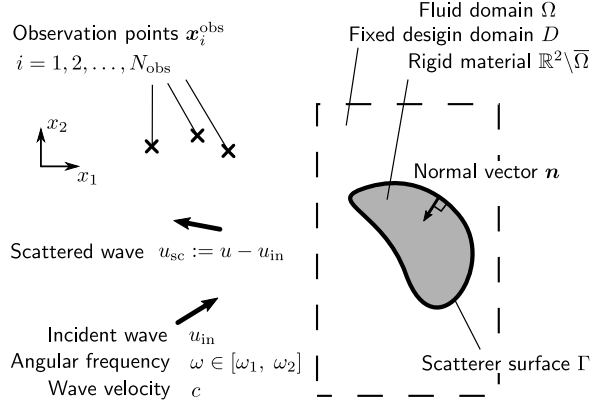}}
 \caption{Two-dimensional acoustic scattering problem.}
 \label{configuration}
\end{figure}

We consider a two-dimensional acoustic scattering caused by acoustically rigid material as shown in \figref{configuration}. We assume that the sound pressure $u^*$ oscillates time-harmonically as $u^*(\bs{x},t)=\Re\left[u(\bs{x})\mathrm{e}^{-\mathrm{i}\omega t}\right]$ where $\bs{x}$, $t$, $u$, and $\omega$ denote the position, time, complex sound pressure, and angular frequency, respectively. $\Re$ denotes the real part of a complex number. In \figref{configuration}, $u_\mathrm{in}(\bs{x})$ is the incident wave with which the scattered wave is defined as $u_\mathrm{sc}(\bs{x}):=u(\bs{x})-u_\mathrm{in}(\bs{x})$. $u(\bs{x})$ satisfies the two-dimensional Helmholtz equation \eqref{helmholtz_eq} and the Neumann boundary condition \equref{neumann_bc} on the scatterer surface $\Gamma$, and the scattered wave $u_\mathrm{sc}(\bs{x})$ satisfies the outgoing radiation condition \equref{outgoing_rc} as follows:
\begin{align}
 \label{helmholtz_eq}
 \nabla^2u(\bs{x})+\frac{\omega^2}{c^2}u(\bs{x})=0&\quad\text{in}\quad\Omega,\\
 \label{neumann_bc}
 \frac{\partial u}{\partial n}(\bs{x})=0&\quad\text{on}\quad\Gamma,\\
 \label{outgoing_rc}
 \frac{\partial u_\mathrm{sc}}{\partial |\bs{x}|}(\bs{x})-\frac{\mathrm{i}\omega}{c}u_\mathrm{sc}(\bs{x})=o(|\bs{x}|^{-\frac{1}{2}})&\quad\text{as}\quad |\bs{x}|\rightarrow\infty,
\end{align}
where $c$ and $\frac{\partial}{\partial n}$ denote the wave velocity in $\Omega$ and the normal direction derivative, respectively. The normal vector $\bs{n}$ on $\Gamma$ is directed inward to the rigid material.

Next, we define the objective function $J$ in the structural optimisation problem. $N_\mathrm{obs}$ observation points $\bs{x}_\mathrm{obs}^i$ are introduced in the fluid domain $\Omega$, with which the objective function $J$ is defined as follows:
\begin{align}
 \label{objective_function}
 J:=\frac{1}{N_\mathrm{obs}(\omega_2-\omega_1)}\int_{\omega_1}^{\omega_2}f(\omega)\mathrm{d}\omega,&\\
 \label{integrand_function}
 f(\omega):=\frac{1}{2}\sum_{i=1}^{N_{\mathrm{obs}}}|u(\omega; \bs{x}_\mathrm{obs}^i)|^2,&
\end{align}
where $[\omega_1,\omega_2]$ is the target band in which the optimised acoustic device will be used. $J$ in \equref{objective_function} physically indicates the sound intensity averaged over the frequency band $[\omega_1,\omega_2]$ and the observation points $\bs{x}_\mathrm{obs}^1,\ldots,\bs{x}_\mathrm{obs}^{N_\mathrm{obs}}$, where $\omega_1$, $\omega_2$, and $\bs{x}_\mathrm{obs}^i\;(i=1,\ldots,N_\mathrm{obs})$ are given by users according to their purpose.
In \equref{integrand_function}, we explicitly list $\omega$ as the argument of $u$ to emphasise that $u$ is a function of the angular frequency. To solve the structural optimisation problem optimising $J$ in \equref{objective_function} under the constraints Eqs. \eqref{helmholtz_eq}--\eqref{outgoing_rc} and $\mathbb{R}^2\backslash\overline{\Omega}\subset D$, where $D\subset\mathbb{R}^2$ is a bounded domain called fixed design domain, we use a topology optimisation method with the level-sets of a B-spline surface~\cite{isakari2017Bspline}. To make the paper self-contained, we give a brief description of the method in \ref{level_set_section}.

\section{Evaluation of the objective function}
  \label{objective_func_section}
In this section, we describe numerical methods to evaluate the objective function $J$ in \equref{objective_function}. Since the integrand $f(\omega)$ cannot generally be written explicitly as a function of $\omega$, the approximate function of $f(\omega)$ is used instead of $f(\omega)$ itself. Here, we may have two possibilities in the approximation; one is to apply the Pad\'{e} approximation directly to $f(\omega)$, while the other is to $u(\omega;\bs{x}_\mathrm{obs}^i)$ on each observation point $\bs{x}_\mathrm{obs}^i$ with which the approximation of $f(\omega)$ is obtained.

\subsection{Direct approximation method}
In this subsection, we present the method in which the integrand $f(\omega)$ is directly approximated. The angular frequency derivatives of $f(\omega)$ at $\omega_0\in[\omega_1, \omega_2]$, the choice of which is discussed later in Section \ref{range_division}, is computed as:
\begin{equation}
 \label{f_freq_der}
  f^{(n)}(\omega_0)=\frac{1}{2}\sum_{i=1}^{N_\mathrm{obs}}\sum_{j=0}^{n} \binom{n}{j}u^{(j)}(\omega_0;\bs{x}_\mathrm{obs}^i)\overline{u}^{(n-j)}(\omega_0;\bs{x}_\mathrm{obs}^i),
\end{equation}
where $\binom{n}{j}$, $g^{(n)}$, and $\overline{g}$ denote the binomial coefficient, $n$-th angular frequency derivative of $g$, and the complex conjugate of $g$, respectively. Numerical procedures to compute $u^{(n)}$ shall be presented in Section \ref{freq_der_section}. The $n$-th order Taylor's coefficient at $\omega_0$ is computed as $f^{(n)}(\omega_0)/n!$ with which the $[M,N]$ Pad\'{e} approximation of $f(\omega)$ at $\omega_0$ is obtained as
\begin{equation}
  \label{fpade}
  f_{[M,N]}(\omega)=\frac{\displaystyle{\sum_{i=0}^M} p_i(\omega-\omega_0)^i}{\displaystyle{\sum_{i=0}^N} q_i(\omega-\omega_0)^i},\quad p_i\;\text{and}\;q_i\in\mathbb{R}.
\end{equation}
See \ref{pade_section} for the definition of the coefficients $p_i$ and $q_i$. We can obtain the approximation of $J$ by integrating \equref{fpade} over the target band $[\omega_1,\omega_2]$. The partial fraction decomposition simplifies this procedure. First, in the case of $M\ge N$, \equref{fpade} is transformed to
\begin{equation}
  \label{poly_division}
  \begin{aligned}
   &f_{[M,N]}(\omega)=\sum_{i=0}^{M-N}r_i(\omega-\omega_0)^i+\frac{\displaystyle{\sum_{i=0}^{N-1}}p'_i(\omega-\omega_0)^i}{\displaystyle{\sum_{i=0}^{N}}q_i(\omega-\omega_0)^i},
  \end{aligned}
\end{equation}
by the polynomial division, where $r_i$ and $p'_i\in\mathbb{R}$ are obtained by Algorithm \ref{division_alg}. In the case of $M<N$, this transformation is not necessary.
\begin{algorithm}
  \caption{The division of polynomials}
  \label{division_alg}
  \begin{algorithmic}[1]
    \REQUIRE $M\ge N$, $q_{N}\ne 0$
    \FOR{$i=M-N$ to $0$}
      \STATE $r_i\leftarrow p_{i+N}/q_{N}$
      \FOR{$j=0$ to $N$}
        \STATE $p_{i+j}\leftarrow p_{i+j}-r_iq_j$
      \ENDFOR
    \ENDFOR
    \FOR{$i=0$ to $N-1$}
      \STATE $p'_i\leftarrow p_i$
    \ENDFOR
  \end{algorithmic}
\end{algorithm}
The partial fraction decomposition is then performed on the second term of the RHS of \equref{poly_division} to have
\begin{equation}
 \label{fpade_decomp}
  \begin{aligned}
   &f_{[M,N]}(\omega)=\sum_{i=0}^{M-N}r_i(\omega-\omega_0)^i+\sum_{i=1}^N \frac{A_i}{\omega-\alpha_i},
  \end{aligned}
\end{equation}
with $\alpha_i$ and $A_i\in\mathbb{C}$.
In \equref{fpade_decomp}, all the poles $\alpha_i$ of \equref{fpade} are required, which are computed by the Durand-Kerner-Aberth (DKA) method. This is a numerical procedure to compute all the zeros of a polynomial at once and thus can be used to find the poles since they are nothing but the zeros of the denominator of \equref{fpade}. In the actual computation, we used an accelerated version of the DKA method with third-order convergence~\cite{aberth1973dka} guaranteed. The coefficients $A_i\;(i=1,\ldots,N)$ of the partial fractions are computed as
\begin{equation}
 \label{heaviside_cover_up}
  A_i=\frac{\displaystyle{\sum_{j=0}^{N-1}}p'_j(\alpha_i-\omega_0)^j}{q_{N}\displaystyle{\prod_{\substack{j=1\\j\ne i}}^N}(\alpha_i-\alpha_j)},
\end{equation}
by the Heaviside cover-up method. Here, it is assumed that there are no multiple roots. We may not guarantee that this assumption always holds, but have not experienced any problems caused by this assumption so far. Then, \equref{fpade_decomp} is easily integrated as
\begin{align}
  \label{J_direct}
    J\simeq\frac{1}{N_\mathrm{obs}(\omega_2-\omega_1)}&\left[\sum_{i=0}^{M-N}\frac{r_i}{i+1}\displaystyle{(\omega-\omega_0)^{i+1}}+\sum_{i=1}^N A_i\displaystyle{\mathrm{Log}(\omega-\alpha_i)}\right]_{\omega_1}^{\omega_2}.
\end{align}
Note that the RHS contains the complex numbers $A_i$ and $\alpha_i$, while the integrand $f(\omega)$ is a real function. It seems strange to introduce complex numbers to evaluate the real function, but this considerably simplifies the procedure to evaluate the objective function.

\subsection{Indirect approximation method}
  \label{indirect_approx_method}
In the other approach, the Pad\'{e} approximation is applied to the frequency response of the complex sound pressure $u(\omega;\bs{x}_\mathrm{obs}^i)$, and then $f(\omega)$ is approximated by the result. For simplicity, we show the case in which one observation point $\bs{x}_\mathrm{obs}$ defines the objective function $J$. We can formulate the indirect approximation method in the same manner for the case that the objective function consists of the sound pressures at multiple observation points. The $[M,N]$ Pad\'{e} approximation of $u(\omega;\bs{x}_\mathrm{obs})$ is obtained as
\begin{equation}
  \label{upade}
  u_{[M,N]}(\omega;\bs{x}_\mathrm{obs})=\frac{\displaystyle{\sum_{i=0}^M} p_i(\omega-\omega_0)^i}{\displaystyle{\sum_{i=0}^N} q_i(\omega-\omega_0)^i},
\end{equation}
with $p_i$ and $q_i\in\mathbb{C}$. $f(\omega)$ in \equref{integrand_function} is then approximated with \equref{upade} as follows:
\begin{align}
  \label{fpade2_before_expanding}
  f_{[2M,2N]}(\omega)&=\frac{1}{2}\frac{\left[\displaystyle{\sum_{i=0}^M} p_i(\omega-\omega_0)^i\right] \left[\displaystyle{\sum_{i=0}^M} \overline{p_i}(\omega-\omega_0)^i\right]}{\left[\displaystyle{\sum_{i=0}^N} q_i(\omega-\omega_0)^i\right] \left[\displaystyle{\sum_{i=0}^N} \overline{q_i}(\omega-\omega_0)^i\right]}\\
  \label{fpade2}
  &=\frac{\displaystyle{\sum_{i=0}^{2M}} \hat{p}_i(\omega-\omega_0)^i}{\displaystyle{\sum_{i=0}^{2N}} \hat{q}_i(\omega-\omega_0)^i},
\end{align}
with $\hat{p}_i$ and $\hat{q}_i\in\mathbb{R}$. In \equref{fpade2_before_expanding}, we used the notation $f_{[2M,2N]}$ to express a rational polynomial approximation of $f$ while $f_{[2M,2N]}$ does not follow the definition of the Pad\'{e} approximation (See \ref{pade_section}). In deriving \equref{fpade2_before_expanding}, we used the fact that $\omega$ and $\omega_0$ are real numbers. $\hat{p}_i$ and $\hat{q}_i$ are obtained by expanding the numerator and denominator of \equref{fpade2_before_expanding} as
\begin{align}
  \label{hat_p}
  \hat{p}_i&=\frac{1}{2}\sum_{j=\max(0,i-M)}^{\min(i,M)} p_j\overline{p_{i-j}}\quad (0\le i\le 2M),\\
  \label{hat_q}
  \hat{q}_i&=\sum_{j=\max(0,i-N)}^{\min(i,N)} q_j\overline{q_{i-j}}\quad (0\le i\le 2N).
\end{align}
We can analytically integrate \equref{fpade2} in a similar manner to the direct approximation method. The objective function $J$ is finally written as:
\begin{align}
  \label{J_indirect}
  J\simeq\frac{1}{\omega_2-\omega_1}&\left[ \sum_{i=0}^{2M-2N}\frac{r_i}{i+1}\displaystyle{(\omega-\omega_0)^{i+1}} +\sum_{i=1}^{2N} A_i\displaystyle{\mathrm{Log}(\omega-\alpha_i)}\right]_{\omega_1}^{\omega_2}
\end{align}
where $r_i$, $A_i$, and $\alpha_i$ are obtained by the partial fraction decomposition of \equref{fpade2}. Equation \eqref{fpade2} has $2N$ poles $\alpha_i$ half of which are the poles of \equref{upade}, and the others are the complex conjugates of them. To see this, let $\alpha_i\;(i=1,\ldots,N)$ be the poles of \equref{upade}; i.e. we have
\begin{equation}
  \label{pole_upade}
  \sum_{i=0}^{N}q_i(\omega-\omega_0)^i=q_N\prod_{i=1}^{N}(\omega-\alpha_i).
\end{equation}
When $\omega$ and $\omega_0$ are real, the complex conjugate of \equref{pole_upade} is
\begin{equation}
  \sum_{i=0}^{N}\overline{q_i}(\omega-\omega_0)^i=\overline{q_N}\prod_{i=1}^{N}(\omega-\overline{\alpha_i}).
\end{equation}
$\overline{\alpha_i}\;(i=1,\ldots,N)$ are, therefore, the poles of \equref{fpade2} as well as $\alpha_i$. Thus, we first seek the poles of \equref{upade} by the DKA method and then compute the complex conjugates of them for the sake of efficiency.

In the case where the multiple observation points exist, we compute the Pad\'{e} approximation \equref{upade} for all the observation points $\bs{x}_\mathrm{obs}^i\;(i=1,2,\ldots,N_\mathrm{obs})$ and take summation for all $\bs{x}_\mathrm{obs}^i$ in Eqs. \eqref{fpade2} and \eqref{J_indirect}.

\section{Comparison of the two methods}\label{sec:comparison_of_the_two_methods}
\begin{figure}
  \centering
  \scalebox{0.75}{\includegraphics{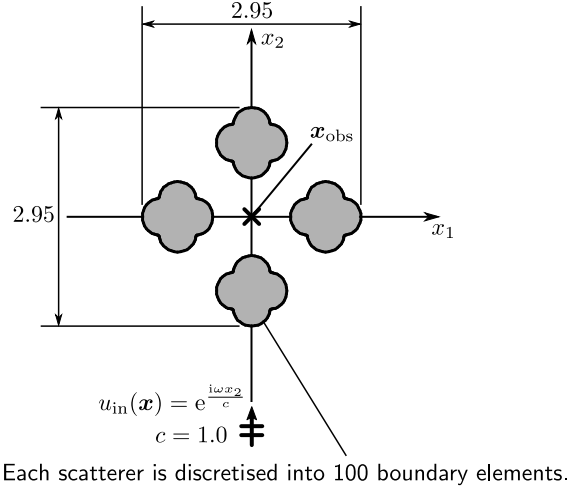}}
  \caption{The observation point surrounded by four rigid scatterers.}
  \label{trochoid4}
\end{figure}
\begin{figure}
  \centering
  \scalebox{0.4}{\includegraphics{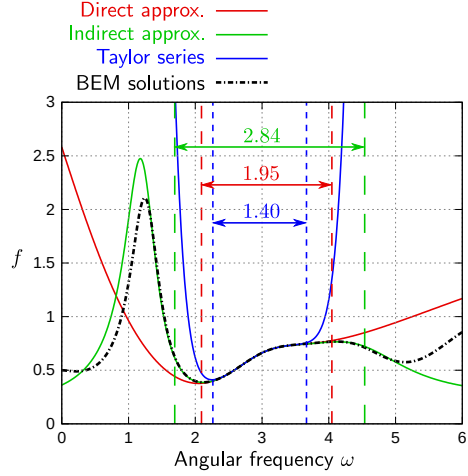}}
  \caption{Sound intensity at $\bs{x}_\mathrm{obs}$ in \figref{trochoid4} estimated by the present approximation methods. $[4,4]$ Pad\'{e} approximation is used in ``Direct approx.'' and ``Indirect approx.''. ``Taylor series'' is truncated at the 8th order. A pair of the same coloured dashed lines show the range where the absolute error of the corresponding approximation to ``BEM solutions'' is less than 0.01.}
  \label{f_plot}
\end{figure}
\begin{figure}
  \centering
  \scalebox{0.4}{\includegraphics{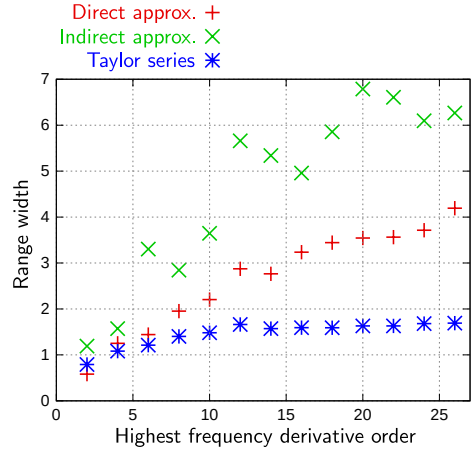}}
  \caption{The horizontal axis shows the highest order of the angular frequency derivative used in each approximation, and the vertical axis shows the width of the frequency range around $\omega_0$ where the absolute error in approximation is less than 0.01.}
  \label{approx_range}
\end{figure}
In this section, we demonstrate which of the direct or indirect approximation methods are better in approximating the frequency-dependent objective function $f(\omega)$. Let us consider the two-dimensional scattering problem as shown in \figref{trochoid4}. Four rigid scatterers, each of which is discretised into 100 boundary elements, are arranged around the origin at which the observation point is placed. The incident wave is the plane wave with the wave velocity $c=1.0$ and unit amplitude propagating along the $x_2$-axis. Figure \ref{f_plot} shows the approximation curves of $f(\omega)$ obtained by the direct and indirect approximation methods with $[4,4]$ Pad\'{e} approximation at $\omega_0=3.0$. For reference, the values of $f(\omega)$ evaluated by the BEM for various $\omega$ in the range (``BEM solutions'') and the Taylor series of $f(\omega)$ at $\omega_0=3.0$ truncated at the 8th order (``Taylor series'') are also plotted in \figref{f_plot}. A pair of the same coloured dashed vertical lines show the range where the absolute error in the approximation to ``BEM solutions'' is less than 0.01. We can see that the indirect approximation method can approximate $f(\omega)$ in the widest region. The direct approximation method is superior to the Taylor expansion while inferior to the indirect approximation method. Figure \ref{approx_range} shows the relationship between the highest order of the angular frequency derivative and the range width where the absolute error is less than 0.01. The order of the Pad\'{e} approximation is chosen such that $M$ and $N$ are identical to each other. The indirect approximation method most rapidly converges to $f(\omega)$. Based on Figs. \ref{f_plot} and \ref{approx_range}, we conclude that the indirect approximation method is better in approximating $f(\omega)$ and henceforth use this method.

\section{Angular frequency derivatives of the solution of the Helmholtz equation}
  \label{freq_der_section}
The angular frequency derivatives $u^{(n)}:=\mathrm{d}^nu/\mathrm{d}\omega^n$ are required to obtain the Pad\'{e} approximation of the angular frequency response. We have already established a low-cost method to compute $u^{(n)}$~\cite{qin2021robust}, which is briefly reviewed in this section. Equations \eqref{helmholtz_eq}--\eqref{outgoing_rc} with the angular frequency $\omega_0$ are equivalent to the following boundary integral equation (BIE):
\begin{equation}
  \label{normal_bie}
  \frac{1}{2}u(\bs{x})+\int_{\Gamma} \frac{\partial G}{\partial n_y}(\bs{x},\bs{y})u(\bs{y}) \mathrm{d}\Gamma_y=u_\mathrm{in}(\bs{x})\quad\bs{x}\in\Gamma,
\end{equation}
where $\frac{\partial}{\partial n_y}$ and $\int_\Gamma\cdot\:\mathrm{d}\Gamma_y$ denote the normal derivative and boundary integral on $\Gamma$ with respect to $\bs{y}$, respectively. Also, $G(\bs{x},\bs{y})$ is the following Green function of the two-dimensional Helmholtz equation:
\begin{equation}
  \label{green_func}
  G(\bs{x},\bs{y}):=\frac{\mathrm{i}}{4}H^{(1)}_0\left(\frac{\omega_0}{c}|\bs{x}-\bs{y}|\right),
\end{equation}
that satisfies the outgoing radiation condition \equref{outgoing_rc}, where $H^{(1)}_0$ is the Hankel function of the first kind and zeroth order. The piecewise constant elements and collocation method, for example, give one numerical solution of \equref{normal_bie}. The solution of \equref{normal_bie} may, however, be polluted by the so-called fictitious eigenfrequency~\cite{burton1971BMmethod}. An alternative BIE free from it is well-established as
\begin{equation}
  \label{burton_miller_bie}
  \begin{aligned}
   \frac{1}{2}u(\bs{x})+\mathrm{p.f.}\int_{\Gamma} W(\bs{x},\bs{y})u(\bs{y}) \mathrm{d}\Gamma_y
   =u_\mathrm{in}(\bs{x})+\gamma q_\mathrm{in}(\bs{x})\quad\bs{x}\in\Gamma,
  \end{aligned}
\end{equation}
where p.f. indicates the finite part of the diverging integral, and the kernel function $W$ is defined as
\begin{equation}
  W(\bs{x},\bs{y}):=\frac{\partial G}{\partial n_y}(\bs{x},\bs{y})+\gamma \frac{\partial^2 G}{\partial n_x\partial n_y}(\bs{x},\bs{y}),
\end{equation}
where $\gamma$ is a complex constant whose imaginary part is non-zero and is set as $\gamma=-\mathrm{i}c/\omega_0$. $q_\mathrm{in}$ denotes the normal derivative of $u_\mathrm{in}$. The BIE \eqref{burton_miller_bie} is called the Burton-Miller type. By differentiating \equref{burton_miller_bie} with respect to angular frequency $n$ times, we obtain the following boundary integral equation:
\begin{equation}
  \label{freq_der_bie}
  \begin{aligned}
    &\frac{1}{2}u^{(n)}(\bs{x})+\mathrm{p.f.}\int_{\Gamma} W(\bs{x},\bs{y})u^{(n)}(\bs{y}) \mathrm{d}\Gamma_y\\
    &=u_\mathrm{in}^{(n)}(\bs{x})+\gamma q_\mathrm{in}^{(n)}(\bs{x})
   -\sum_{i=0}^{n-1}\binom{n}{i} \mathrm{p.f.}\int_{\Gamma} W^{(n-i)}(\bs{x},\bs{y})u^{(i)}(\bs{y}) \mathrm{d}\Gamma_y.
  \end{aligned}
\end{equation}
We can obtain $u^{(0)}(=u)$ by solving \equref{burton_miller_bie} and recursively solve \equref{freq_der_bie} from $n=1$ to the highest ($=(M+N)$-th in this paper) derivative in order. By solving Eqs.~\eqref{burton_miller_bie} and \eqref{freq_der_bie} at $\omega=\omega_0$, we can compute the boundary values of $u(\bs{x};\omega_0)$ and $u^{(n)}(\bs{x};\omega_0)$. The forward mode automatic differentiation \cite{qin2021robust} is available to compute the higher-order angular frequency derivatives of $W$, $u_\mathrm{in}$ and $q_\mathrm{in}$ in \equref{freq_der_bie}. All the linear equations obtained by the discretised \equref{freq_der_bie} share the identical LHS matrix for all $n$. It is, therefore, reasonable to use a direct solver, such as the LU factorisation, to solve them. We use the hierarchical matrix method~\cite{mario2008hmatrix} to efficiently assemble the LHS and LU factorise it.

The objective function $J$ and its sensitivity (to be discussed in Section \ref{sensitivity_section}) are evaluated via $u(\bs{x})$ and $u^{(n)}(\bs{x})$ $\bs{x}\in\Omega$, which is obtained by the following integral representation:
\begin{equation}
  \label{interior_domain_value}
  u(\bs{x})=u_\mathrm{in}(\bs{x})-\int_{\Gamma}\frac{\partial G}{\partial n_y}(\bs{x},\bs{y})u(\bs{y})\mathrm{d}\Gamma_y,\quad\bs{x}\in\Omega
\end{equation}
and its $n$-th angular frequency derivatives. The automatic differentiation is again engaged to evaluate them.

\section{A fast method to compute the RHS of \equref{freq_der_bie}}\label{a_fast_method_to_compute_the_rhs_of}
The RHS of \equref{freq_der_bie} involves the layer potential whose kernels are the angular frequency derivatives of the Green function. It is quite time-consuming to compute them naively with the automatic differentiation. We can accelerate such computations by the fast multipole method (FMM) \cite{rokhlin1985fmm, greengard1987fmm} with minor modifications. Note that this  also applies to the angular frequency derivatives of \equref{interior_domain_value}. We first introduce an auxiliary frequency-dependent variable $v$ defined as
\begin{equation}
  v^{(i)}(\bs{x})=
  \begin{cases}
    u^{(i)}(\bs{x}) & 0\le i < n\\
    0 & i\ge n
  \end{cases}
  \quad.
\end{equation}
The boundary integrals in the RHS of \equref{freq_der_bie} is then rewritten as follows:
\begin{equation}
  \label{boundary_integrals}
  \begin{aligned}
   \sum_{i=0}^{n-1}\binom{n}{i}\mathrm{p.f.}\int_{\Gamma} W^{(n-i)}(\bs{x},\bs{y})u^{(i)}(\bs{y}) \mathrm{d}\Gamma_y
   =\left[ \mathrm{p.f.}\int_{\Gamma} W(\bs{x},\bs{y})v(\bs{y}) \mathrm{d}\Gamma_y \right]^{(n)}.
  \end{aligned}
\end{equation}
Since the representation in the parentheses of \equref{boundary_integrals} is nothing but a usual layer potential (see \equref{burton_miller_bie}), we can use the standard FMM to compute it. Furthermore, the automatic differentiation provides its angular frequency derivatives simultaneously. Thus, we just need to install the automatic differentiation routine in an existing FMM code. We here summarise the formulas of the low-frequency FMM used in this study.\\
{\bf The multipole moment}:
\begin{equation}
  \label{multipole_moment}
  M_m(\bs{y}_0)=(-1)^m\int_{\Gamma_s}\frac{\partial I_m}{\partial n_y}(\bs{y}-\bs{y}_0)v(\bs{y})\mathrm{d}\Gamma_y,
\end{equation}
in which $\Gamma_s$ is a part of the boundary $\Gamma$.\\
{\bf M2M formula}:
\begin{equation}
  \label{M2M}
  M_m(\bs{y}'_0)=\sum_{k}I_{m-k}(\bs{y}'_0-\bs{y}_0)M_k(\bs{y}_0),
\end{equation}
{\bf M2L formula}:
\begin{equation}
  \label{M2L}
  L_k(\bs{x}_0)=\sum_{m} O_{k-m}(\bs{x}_0-\bs{y}_0)M_m(\bs{y}_0),
\end{equation}
{\bf L2L formula}:
\begin{equation}
  \label{L2L}
  L_m(\bs{x}'_0)=\sum_{k}I_{m-k}(\bs{x}'_0-\bs{x}_0)L_k(\bs{x}_0),
\end{equation}
{\bf The local expansion}:
\begin{equation}
  \label{local_expansion}
  \begin{aligned}
   \int_{\Gamma_s}W(\bs{x},\bs{y})v(\bs{y})\mathrm{d}\Gamma_y
   =\frac{\mathrm{i}}{4}\sum_{k}\left[I_k(\bs{x}-\bs{x}_0)+\gamma \frac{\partial I_k}{\partial n}(\bs{x}-\bs{x}_0)\right]L_{-k}(\bs{x}_0).
  \end{aligned}
\end{equation}
In the above formulas, $I_n$ and $O_n$ are the solutions of the Helmholtz equation in $\mathbb{R}^2$ whose explicit representations are given as
\begin{align}
  \label{I_function}
  I_n(\bs{x})&=\mathrm{i}^nJ_n\left(\frac{\omega_0 r}{c}\right)\mathrm{e}^{\mathrm{i}n\theta},\\
  \label{O_function}
  O_n(\bs{x})&=\mathrm{i}^nH^{(1)}_n\left(\frac{\omega_0 r}{c}\right)\mathrm{e}^{\mathrm{i}n\theta},
\end{align}
where $(r,\theta)$ is the polar coordinate of $\bs{x}$, and $J_n$ is the Bessel function of $n$-th order. Note that $I_n$, $O_n$, $M_n$, and $L_n$ are frequency-dependent and thus their angular frequency derivatives must be computed. They are, however, computed ``automatically'' with the help of the automatic differentiation. In Eqs. \eqref{M2M}--\eqref{L2L}, the infinite sums are truncated at a finite order. We have numerically confirmed that we can determine  the truncation order in a standard way~\cite{coifman1993fast} with the standard FMM for the two-dimensional Helmholtz equation.

\section{Sensitivity analysis}
  \label{sensitivity_section}
In this section, we discuss the computation of the design sensitivity. We adopt the topological derivative of the objective function as the design sensitivity. In this paper, $[\mathcal{D}_\mathrm{T} g](\bs{x})$ denotes the topological derivative of $g$. Let us assume that a small rigid circular scatterer of radius $\varepsilon>0$ is introduced to $\bs{x}\in\Omega$, and a variable $g$ changes to
$g+\delta g(\bs{x};\varepsilon)$
due to the topological perturbation.
In this case, the topological derivative $\mathcal{D}_\mathrm{T}g$ is defined by the following limit:
\begin{equation}
 [\mathcal{D}_\mathrm{T}g](\bs{x}):=\lim_{\varepsilon\downarrow 0}\frac{\delta g(\bs{x};\varepsilon)}{\pi \varepsilon^2}\quad \bs{x}\in\Omega.
\label{def:topological_derivative}
\end{equation}
In our method, the objective function $J$ is computed by the indirect approximation method as \equref{J_indirect}. Its sensitivity is also computed based on this approximation. In this section, we show the formulation in the case of $N_\mathrm{obs}=1$ and discuss the case of $N_\mathrm{obs}>1$ later. The procedure to derive $\mathcal{D}_\mathrm{T}J$ is as follows: first, the topological derivatives $[\mathcal{D}_\mathrm{T}u_\mathrm{obs}](\bs{x})$ and $[\mathcal{D}_\mathrm{T}u_\mathrm{obs}^{(n)}](\bs{x})$, where $u_\mathrm{obs}:=u(\bs{x}_\mathrm{obs},\omega_0)$ is introduced, are computed. Then, $\mathcal{D}_\mathrm{T}J$ shall accordingly be derived by the chain rule.

\subsection{$\mathcal{D}_\mathrm{T}u_\mathrm{obs}$ and $\mathcal{D}_\mathrm{T}u_\mathrm{obs}^{(n)}$}
First, let us define the adjoint variable $\tilde{u}$ by the following boundary value problem:
\begin{align}
  \label{helmholtz_eq_adj}
  \nabla^2 \tilde{u}(\bs{x})+\frac{\omega^2}{c^2}\tilde{u}(\bs{x})+\delta(\bs{x}-\bs{x}_\mathrm{obs})=0\quad&\mathrm{in}\quad\Omega,\\
  \label{neumann_bc_adj}
  \frac{\partial \tilde{u}}{\partial n}(\bs{x})=0\quad&\mathrm{on}\quad\Gamma,\\
  \label{outgoing_rc_adj}
  \frac{\partial \tilde{u}}{\partial |\bs{x}|}(\bs{x})-\frac{\mathrm{i}\omega}{c}\tilde{u}(\bs{x})=o\left(|\bs{x}|^{-\frac{1}{2}}\right)\quad&\mathrm{as}\quad|\bs{x}|\rightarrow\infty,
\end{align}
where $\delta$ denotes the Dirac delta. Equations \eqref{helmholtz_eq_adj}--\eqref{outgoing_rc_adj} are equivalent to the following BIE:
\begin{equation}
  \label{adjoint_bie}
   \begin{aligned}
    \frac{1}{2}\tilde{u}(\bs{x})+\mathrm{p.f.}\int_\Gamma W(\bs{x},\bs{y})\tilde{u}(\bs{y})\mathrm{d}\Gamma_y
    =G(\bs{x},\bs{x}_\mathrm{obs})+\gamma\frac{\partial G}{\partial n}(\bs{x},\bs{x}_\mathrm{obs})\quad \bs{x}\in\Gamma,
   \end{aligned}
\end{equation}
where the Burton-Miller type formulation is adopted. The angular frequency derivatives $\tilde{u}^{(n)}$ are evaluated in the same way as the primal problem (see Section \ref{freq_der_section}). The discretised BIEs share the same LHS matrix as that of \equref{burton_miller_bie}, and thus the LU factorised matrix in the primal problem can be recycled.
$[\mathcal{D}_\mathrm{T}u_\mathrm{obs}](\bs{x})$ is evaluated with $\tilde{u}$ as follows~\cite{carpio2008inverse}:
\begin{equation}
 \label{Du_eval}
  [\mathcal{D}_\mathrm{T}u_\mathrm{obs}](\bs{x})=2\sum_{j=1}^2
  \frac{\partial \tilde{u}(\bs{x})}{\partial x_j}
  \frac{\partial u(\bs{x})}{\partial x_j}
  -\frac{\omega^2}{c^2}\tilde{u}(\bs{x})u(\bs{x}).
\end{equation}
We can obtain the topological derivatives $[\mathcal{D}_\mathrm{T}u_\mathrm{obs}^{(n)}](\bs{x})$ by differentiating \equref{Du_eval} $n$ times with respect to $\omega$. In practice, however, we do not need to differentiate \equref{Du_eval} explicitly because the automatic differentiation is available. Here, we implicitly swapped the topological derivative and angular frequency derivative as:
\begin{equation}
  \label{derivaticve_swap}
  \mathcal{D}_\mathrm{T}[u_\mathrm{obs}^{(n)}]=[\mathcal{D}_\mathrm{T}u_\mathrm{obs}]^{(n)}.
\end{equation}
It has not yet been proved rigorously but inferred numerically that \equref{derivaticve_swap} holds for arbitrary $n$~\cite{qin2021robust}.

\subsection{Topological derivative of the objective function}
In this subsection, we derive $\mathcal{D}_\mathrm{T}J$ based on the chain rule. Let us first recall the procedure to compute  $J$ in the indirect approximation method:
\begin{enumerate}
  \item Obtain the Pad\'{e} approximation of $u(\omega;\bs{x}_\mathrm{obs})$ as \equref{upade},
  \item Obtain the approximation of $f(\omega)$ as \equref{fpade2},
  \item Obtain the partial fraction decomposition of the approximation,
  \item Integrate the approximation and compute $J$ as \equref{J_indirect}.
\end{enumerate}
Similarly, the procedure to compute $\mathcal{D}_\mathrm{T}J$ is as follows:
\begin{enumerate}
  \item Compute the topological derivatives of $p_i$ and $q_i$ in \equref{upade},
  \item Compute the topological derivatives of $\hat{p}_i$ and $\hat{q}_i$ in \equref{fpade2},
  \item Compute the topological derivatives of $r_i$, $A_i$, and $\alpha_i$ in \equref{J_indirect},
  \item Compute the topological derivative of $J$.
\end{enumerate}

\subsubsection{$\mathcal{D}_\mathrm{T}p_i$ and $\mathcal{D}_\mathrm{T}q_i$} \label{Dp_Dq_subsection}
Let us recall that the coefficients of the Pad\'{e} approximation $p_i$ and $q_i$ are obtained by the following linear equations (See \ref{pade_section}):
\begin{align}
  \label{tmp1}
  &p_0=u_\mathrm{obs},\\
  \label{tmp2}
  &q_0=1,\\
  \label{tmp3}
  &p_i-\sum_{j=1}^i\frac{u_\mathrm{obs}^{(i-j)}}{(i-j)!}q_j=\frac{u_\mathrm{obs}^{(i)}}{i!},
\end{align}
for $1\le i\le M+N$, where $p_i:=0\;(i>M)$ and $q_i:=0\;(i>N)$ are defined. The topological derivatives of Eqs. \eqref{tmp1}--\eqref{tmp3} are given as:
\begin{align}
  &\mathcal{D}_\mathrm{T}p_0=\mathcal{D}_\mathrm{T}u_\mathrm{obs},\\
  &\mathcal{D}_\mathrm{T}q_0=0,\\
  &\mathcal{D}_\mathrm{T}p_i-\sum_{j=1}^i\frac{u_\mathrm{obs}^{(i-j)}}{(i-j)!}\mathcal{D}_\mathrm{T}q_j
  =\frac{\mathcal{D}_\mathrm{T}u_\mathrm{obs}^{(i)}}{i!}+\sum_{j=1}^i\frac{\mathcal{D}_\mathrm{T}u_\mathrm{obs}^{(i-j)}}{(i-j)!}q_j,
\end{align}
which are solved in the same way as Eqs. \eqref{tmp1}--\eqref{tmp3} (See \ref{pade_section}).

\subsubsection{$\mathcal{D}_\mathrm{T}\hat{p}_i$ and $\mathcal{D}_\mathrm{T}\hat{q}_i$}
The topological derivatives of Eqs. \eqref{hat_p} and \eqref{hat_q} are respectively expressed as
\begin{align}
  \label{D_hat_p}
  &\mathcal{D}_\mathrm{T}\hat{p}_i=\frac{1}{2}\sum_{j=\max(0,i-M)}^{\min(i,M)} (p_j\overline{\mathcal{D}_\mathrm{T}p_{i-j}}+\mathcal{D}_\mathrm{T}p_j\overline{p_{i-j}}),
\end{align}
for $0\le i\le 2M$, and
\begin{align}
  \label{D_hat_q}
  &\mathcal{D}_\mathrm{T}\hat{q}_i=\sum_{j=\max(0,i-N)}^{\min(i,N)} (q_j\overline{\mathcal{D}_\mathrm{T}q_{i-j}}+\mathcal{D}_\mathrm{T}q_j\overline{q_{i-j}}),
\end{align}
for $0\le i\le 2N$, where $\mathcal{D}_\mathrm{T}\overline{z}=\overline{\mathcal{D}_\mathrm{T}z}$ for $z\in\mathbb{C}$ is used.

\subsubsection{$\mathcal{D}_\mathrm{T}r_i$}
As the result of the polynomial division, \equref{fpade2} is transformed to
\begin{equation}
  f(\omega)=\sum_{i=0}^{2M-2N}r_i(\omega-\omega_0)^i+\frac{\displaystyle{\sum_{i=0}^{2N-1}}\hat{p}'_i(\omega-\omega_0)^i}{\displaystyle{\sum_{i=0}^{2N}}\hat{q}_i(\omega-\omega_0)^i}.
\end{equation}
The algorithm to compute $r_i$, $\hat{p}'_i$, $\mathcal{D}_\mathrm{T}r_i$, and $\mathcal{D}_\mathrm{T}\hat{p}'_i$ simultaneously is summarised as Algorithm \ref{D_r_D_hat_p_alg}.
\begin{algorithm}
  \caption{The division of polynomials with the topological derivatives}
  \label{D_r_D_hat_p_alg}
  \begin{algorithmic}[1]
    \REQUIRE $M\ge N$, $\hat{q}_{2N}\ne 0$
    \FOR{$i=2M-2N$ to $0$}
      \STATE $r_i\leftarrow \hat{p}_{i+2N}/\hat{q}_{2N}$
      \STATE $\mathcal{D}_\mathrm{T}r_i\leftarrow \left(\mathcal{D}_\mathrm{T}\hat{p}_{i+2N}\hat{q}_{2N}-\hat{p}_{i+2N}\mathcal{D}_\mathrm{T}\hat{q}_{2N}\right)/(\hat{q}_{2N})^2$
      \FOR{$j=0$ to $2N$}
        \STATE $\hat{p}_{i+j}\leftarrow \hat{p}_{i+j}-r_i\hat{q}_j$
        \STATE $\mathcal{D}_\mathrm{T}\hat{p}_{i+j}\leftarrow \mathcal{D}_\mathrm{T}\hat{p}_{i+j}-\mathcal{D}_\mathrm{T}r_i\hat{q}_j-r_i\mathcal{D}_\mathrm{T}\hat{q}_j$
      \ENDFOR
    \ENDFOR
    \FOR{$i=0$ to $2N-1$}
      \STATE $\hat{p}'_i\leftarrow \hat{p}_i$
      \STATE $\mathcal{D}_\mathrm{T}\hat{p}'_i\leftarrow \mathcal{D}_\mathrm{T}\hat{p}_i$
    \ENDFOR
  \end{algorithmic}
\end{algorithm}
The 3rd and 6th rows of Algorithm \ref{D_r_D_hat_p_alg} are obtained by the topological derivatives of the 2nd and 5th rows, respectively.

\subsubsection{$\mathcal{D}_\mathrm{T}\alpha_i$}
The set of the poles $\alpha_i$ of \equref{fpade2} consists of the poles of \equref{upade} and the complex conjugates of them. Let $\alpha_i\;(i=1,\ldots,N)$ denote the poles of \equref{upade}, and $\alpha_i=\overline{\alpha_{i-N}}\;(i=N+1,\ldots,2N)$. The former satisfies
\begin{equation}
  \label{tmp}
  \sum_{j=0}^N q_j(\alpha_i-\omega_0)^j=0\quad (i=1,\ldots,N),
\end{equation}
because $\alpha_i$ is a root of the denominator of \equref{upade}. The topological derivative of \equref{tmp} leads to
\begin{equation}
  \label{D_alpha}
  \mathcal{D}_\mathrm{T}\alpha_i=-\frac{\displaystyle{\sum_{j=0}^N} \mathcal{D}_\mathrm{T}q_j(\alpha_i-\omega_0)^j}{\displaystyle{\sum_{j=1}^N} jq_j(\alpha_i-\omega_0)^{j-1}}.
\end{equation}
The latter is obtained as $\mathcal{D}_\mathrm{T}\alpha_i=\overline{\mathcal{D}_\mathrm{T}\alpha_{i-N}}\; (i=N+1,\ldots,2N)$.

\subsubsection{$\mathcal{D}_\mathrm{T}A_i$}
\label{DA_subsection}
The coefficient $A_i$ in \equref{J_indirect} is evaluated by the Heaviside cover-up method as:
\begin{equation}
  \label{heaviside_cover_up2}
  A_i=\frac{\displaystyle{\sum_{j=0}^{2N-1}}\hat{p}'_j(\alpha_i-\omega_0)^j}{\hat{q}_{2N}\displaystyle{\prod_{\substack{j=1\\j\ne i}}^{2N}}(\alpha_i-\alpha_j)}.
\end{equation}
The topological derivative of \equref{heaviside_cover_up2} leads to $\mathcal{D}_\mathrm{T}A_i$:
\begin{equation}
  \begin{aligned}
    \mathcal{D}_\mathrm{T}A_i=\frac{1}{\hat{q}_{2N}\displaystyle{\prod_{\substack{j=1\\j\neq i}}^{2N}} (\alpha_i-\alpha_j)}
    \left[
      \sum_{j=0}^{2N-1}\mathcal{D}_\mathrm{T}\hat{p}'_j(\alpha_i-\omega_0)^j
   +\mathcal{D}_\mathrm{T}\alpha_i\sum_{j=1}^{2N-1}j\hat{p}'_j(\alpha_i-\omega_0)^{j-1}\right.\\
    \left.
      -\left(
        \frac{\mathcal{D}_\mathrm{T}\hat{q}_{2N}}{\hat{q}_{2N}}
        +\sum_{\substack{j=1\\j\neq i}}^{2N}\frac{\mathcal{D}_\mathrm{T}\alpha_i-\mathcal{D}_\mathrm{T}\alpha_j}{\alpha_i-\alpha_j}
      \right)\sum_{j=0}^{2N-1}\hat{p}'_j(\alpha_i-\omega_0)^j
    \right].
  \end{aligned}
\end{equation}

\subsubsection{$\mathcal{D}_\mathrm{T}J$}
Finally, the topological derivative of \equref{J_indirect} gives the estimation of $\mathcal{D}_\mathrm{T}J$:
\begin{equation}
  \label{DJ_indirect}
  \begin{aligned}
    \mathcal{D}_\mathrm{T}J\simeq&\frac{1}{\omega_2-\omega_1}\left[\sum_{i=0}^{2M-2N}\frac{\mathcal{D}_\mathrm{T}r_i}{i+1}(\omega-\omega_0)^{i+1}\right.\\
              &\left.+\sum_{i=1}^{2N}\mathcal{D}_\mathrm{T}A_i\mathrm{Log}(\omega-\alpha_i)
                      -\sum_{i=1}^{2N}A_i\frac{\mathcal{D}_\mathrm{T}\alpha_i}{\omega-\alpha_i}\right]_{\omega_1}^{\omega_2}.
  \end{aligned}
\end{equation}
Note that we have so far discussed the case of $N_\mathrm{obs}=1$.
In the case of $N_\mathrm{obs}>1$, the topological derivatives $\mathcal{D}_\mathrm{T}p_i$, $\mathcal{D}_\mathrm{T}q_i$, $\mathcal{D}_\mathrm{T}\hat{p}_i$, $\mathcal{D}_\mathrm{T}\hat{q}_i$, $\mathcal{D}_\mathrm{T}r_i$, $\mathcal{D}_\mathrm{T}\alpha_i$, and $\mathcal{D}_\mathrm{T}A_i$ are computed for each observation point, and summations are taken for all the observation points in \equref{DJ_indirect}.

\section{Adaptive subdivision of the target band}
\label{range_division}
In Sections \ref{objective_func_section} and \ref{sensitivity_section}, we have shown a way to estimate the objective function $J$ and its topological derivative $\mathcal{D}_\mathrm{T}J$ with the Pad\'{e} approximation. The range where the approximation holds is, however, limited to a certain range of the frequency around the centre of the approximation $\omega_0$. If the target band $[\omega_1,\omega_2]$ exceeds the limit, $J$ and $\mathcal{D}_\mathrm{T}J$ can no longer be estimated accurately. In this case, the target band is divided into some subintervals, and the approximations are computed in each interval. In this section, we discuss a reasonable way to determine the number of subintervals adaptively in the optimisation.
When determining the number of subintervals, we monitor the approximation accuracy of $\mathcal{D}_\mathrm{T}J$ rather than $J$ itself because the range where $\mathcal{D}_\mathrm{T}J$ is accurately approximated is narrower than the corresponding range of $J$, which will be observed in Section \ref{example_section}.

Let us recall that, in the case of $N_\mathrm{obs}=1$, $\mathcal{D}_\mathrm{T}J$ is represented as
    \begin{equation}
      \mathcal{D}_\mathrm{T}J=\frac{1}{\omega_2-\omega_1}\int_{\omega_1}^{\omega_2}\mathcal{D}_\mathrm{T}f(\omega)\text{d}\omega,
    \end{equation}
    in which $\mathcal{D}_\mathrm{T}f(\omega)$ is approximated as
    \begin{equation}
      \label{Dfpade}
      \begin{aligned}
        \mathcal{D}_\mathrm{T}f_{[2M,2N]}(\omega)&=
          \frac{\displaystyle{\sum_{i=0}^{2M}} \mathcal{D}_\mathrm{T}\hat{p}_i(\omega-\omega_0)^i}{\displaystyle{\sum_{i=0}^{2N}} \hat{q}_i(\omega-\omega_0)^i}\\
          &-\frac{\left[\displaystyle{\sum_{i=0}^{2M}} \hat{p}_i(\omega-\omega_0)^i\right]\left[\displaystyle{\sum_{i=0}^{2N}} \mathcal{D}_\mathrm{T}\hat{q}_i(\omega-\omega_0)^i\right]}{\left[\displaystyle{\sum_{i=0}^{2N}} \hat{q}_i(\omega-\omega_0)^i\right]^2},
      \end{aligned}
    \end{equation}
    by the indirect approximation. Let us here define an angular frequency range $[\omega_\mathrm{L}, \omega_\mathrm{R}]$ in which $\mathcal{D}_\mathrm{T}f_{[2M,2N]}$ gives a sufficiently accurate approximation of $\mathcal{D}_\mathrm{T}f$ as follows:
    \begin{align}
      \label{omega_R_L_def1}
      &\max\;\omega_\mathrm{R}-\omega_\mathrm{L}\\
      &\quad \text{subject\; to}\nonumber\\
      \label{omega_R_L_def_3}
      &\quad e_i(\omega)\le\delta_{\mathcal{D}_\mathrm{T}f}\quad\forall\omega\in[\omega_\mathrm{L},\omega_\mathrm{R}]\;\forall\bs{x}_i,\\
      &\quad[\omega_\mathrm{L},\omega_\mathrm{R}]\subset[\omega_1,\omega_2],\\
      &\quad\omega_0\in[\omega_\mathrm{L},\omega_\mathrm{R}], \\
      &\quad e_i(\omega):=|\mathcal{D}_\mathrm{T}f_{[2M,2N]}(\omega;\bs{x}_i)-\mathcal{D}_\mathrm{T}f_{[2(M-1),2(N-1)]}(\omega;\bs{x}_i)|
      \label{omega_R_L_def2}
    \end{align}
    where $\delta_{\mathcal{D}_\mathrm{T}f}$ and $\bs{x}_i$ denote the tolerance given by the user and the points where the computation of the topological derivative is required by the optimiser, respectively. $\mathcal{D}_\mathrm{T}f_{[2(M-1),2(N-1)]}$ is derived from the $[M-1,N-1]$ Pad\'{e} approximation $u_{[M-1,N-1]}$. We can estimate that the approximation of the topological derivative is valid in $[\omega_\mathrm{L},\omega_\mathrm{R}]$ because the results of the two approximations with different degrees are consistent. To find $[\omega_\mathrm{L},\omega_\mathrm{R}]$, we define $[\omega_\mathrm{L}^i,\omega_\mathrm{R}^i]$ for each $\bs{x}_i$ as follows:
    \begin{align}
        &\max\;\omega_\mathrm{R}^i-\omega_\mathrm{L}^i\\
        &\quad \text{subject\; to}\nonumber\\
        &\quad e_i(\omega)\le\delta_{\mathcal{D}_\mathrm{T}f},\quad\forall\omega\in[\omega_\mathrm{L}^i,\omega_\mathrm{R}^i],\\
        &\quad[\omega_\mathrm{L}^i,\omega_\mathrm{R}^i]\subset[\omega_1,\omega_2],\\
        &\quad\omega_0\in[\omega_\mathrm{L}^i,\omega_\mathrm{R}^i].
    \end{align}
    If $e_i(\omega)$ monotonically increases as $\omega$ goes away from $\omega_0$, which is true when $\omega$ is sufficiently close to $\omega_0$, we can find $\omega_\mathrm{L}^i$ by solving the nonlinear equation $e_i(\omega)=\delta_{\mathcal{D}_\mathrm{T}f}$ in $[\omega_1,\omega_0]$, and so is $\omega_\mathrm{R}^i$. We use the bisection to find them. In the case of $e_i(\omega_1)<\delta_{\mathcal{D}_\mathrm{T}f}$, $\omega_\mathrm{L}^i=\omega_1$. $[\omega_\mathrm{L},\omega_\mathrm{R}]$ is equal to the overlapping region of $[\omega_\mathrm{L}^i,\omega_\mathrm{R}^i]$ for all $\bs{x}_i$.

The number of the subintervals $N_\mathrm{div}$ is determined with $[\omega_\mathrm{L},\omega_\mathrm{R}]$ as follows:
    \begin{align}
      \label{define_N_div1}
      N_\mathrm{div}'&=\left\lceil \frac{\omega_2-\omega_1}{\omega_\mathrm{R}-\omega_\mathrm{L}}\right\rceil,\\
      \label{define_N_div2}
      N_\mathrm{div}&=
      \begin{cases}
        N_\mathrm{div}' & \text{if }N_\mathrm{div}'\text{ is odd},\\
        N_\mathrm{div}'+1 & \text{if }N_\mathrm{div}'\text{ is even},
      \end{cases}
    \end{align}
    where $\lceil\cdot\rceil$ denotes the ceiling function. We solve the primal and adjoint problems at $\omega_0:=\frac{1}{2}(\omega_1+\omega_2)$, find $[\omega_\mathrm{L},\omega_\mathrm{R}]$ defined by Eqs. \eqref{omega_R_L_def1}--\eqref{omega_R_L_def2}, and determine $N_\mathrm{div}$ by \equref{define_N_div2}.
    The target band $[\omega_1,\omega_2]$ is divided into $N_\mathrm{div}$ equally spaced subintervals. In the analyses in the subintervals, the centre of the Pad\'{e} approximation is set to the midpoint of each interval. Note that we can recycle the Pad\'{e} approximation at $\frac{1}{2}(\omega_1+\omega_2)$ because $N_\mathrm{div}$ is always odd owing to its definition \equref{define_N_div2}.

\section{Numerical examples}\label{example_section}
\subsection{Topological derivative}
\begin{figure}
 \centering
 \scalebox{0.85}{\includegraphics{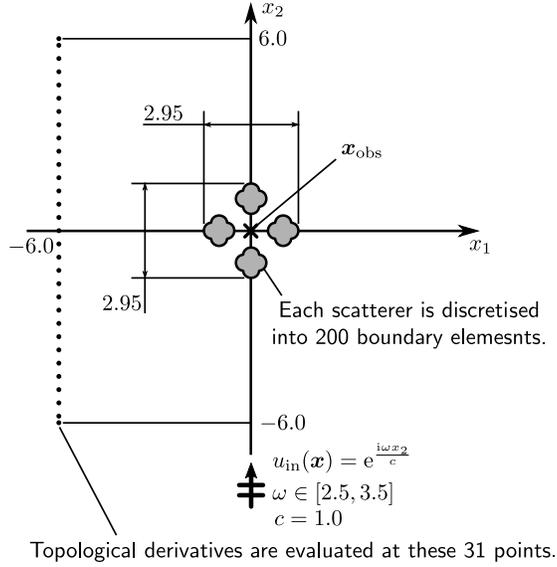}}
 \caption{Scatterers surrounding the observation point and the points where topological derivatives are evaluated.}
 \label{section7subsec1_configuration}
\end{figure}
\begin{figure}
  \centering
  \scalebox{0.4}{\includegraphics{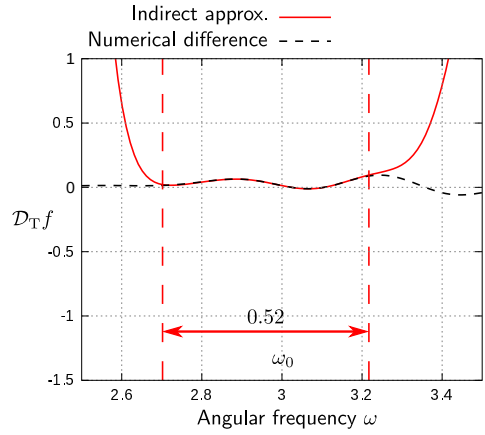}}\\
  (a) $\bs{x}=(-6.0,-6.0)$\\
  \vspace{2mm}
  \scalebox{0.4}{\includegraphics{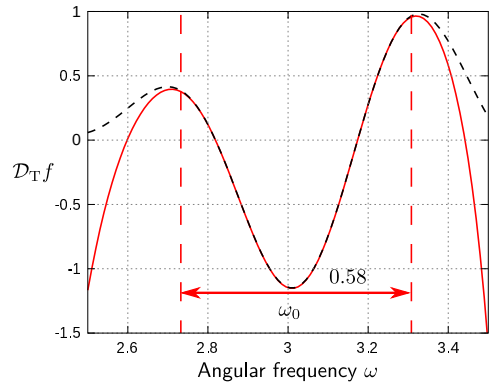}}\\
  (b) $\bs{x}=(-6.0,0.0)$
  \caption{The plots of $[\mathcal{D}_\mathrm{T}f](\bs{x})$ evaluated by \equref{Dfpade} (``Indirect approx.'') and \equref{numerical_difference} (``Numerical difference''). A pair of vertical dashed lines show the range where the absolute error of ``Indirect approx.'' against ``Numerical difference'' is less than 0.01.}
  \label{Df_pt1_466}
\end{figure}
\begin{figure}
  \centering
  \scalebox{0.4}{\includegraphics{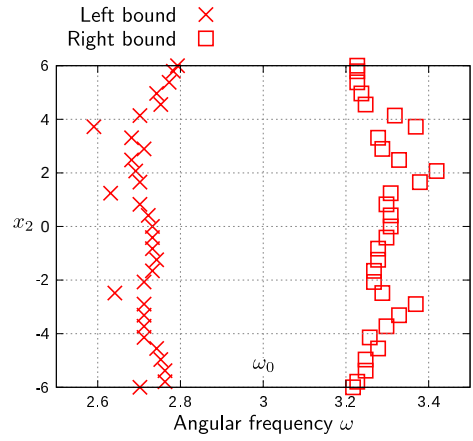}}
  \caption{The range where the absolute error of $\mathcal{D}_\mathrm{T}f$ in ``Indirect approx.'' against ``Numerical difference'' is less than 0.01.}
  \label{range_comparison}
\end{figure}
\begin{figure}
  \centering
  \scalebox{0.4}{\includegraphics{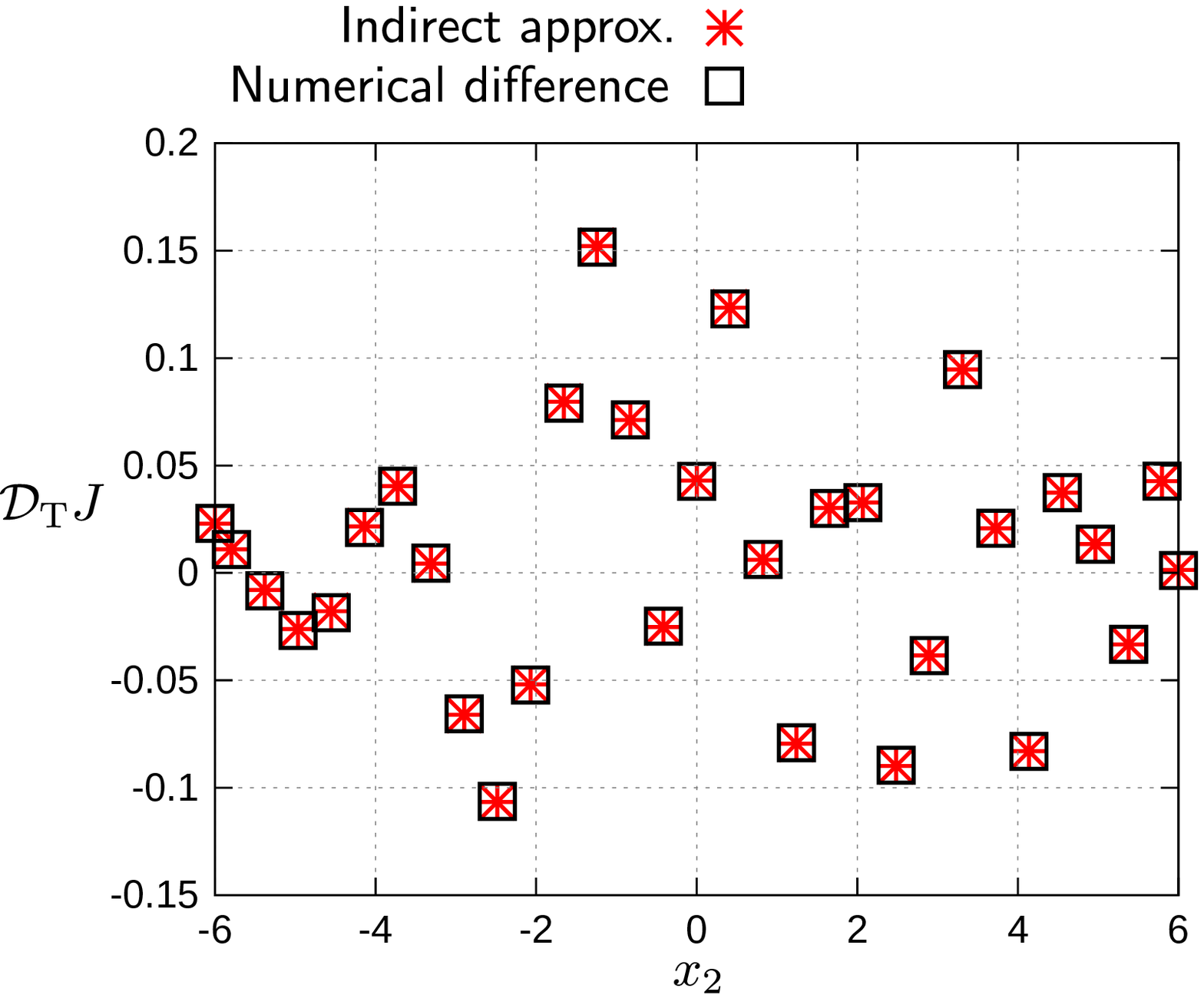}}
  \caption{The topological derivative of the objective function evaluated by \equref{DJ_indirect} (``Indirect approx.'') and \equref{numerical_difference} (``Numerical difference'').}
  \label{DJ_plot}
\end{figure}
In this subsection, we first numerically validate \equref{Dfpade} and check the angular frequency range where the approximation holds. We then validate the topological derivative \equref{DJ_indirect} computed with the target band subdivision described in Section \ref{range_division}.
In the validations, as a reference, we use the following finite numerical difference of the relevant quantity:
\begin{equation}
  \label{numerical_difference}
  [\mathcal{D}_\mathrm{T}g](\bs{x})\simeq \frac{\delta g(\bs{x};\varepsilon)}{\pi\varepsilon^2},
\end{equation}
where $\delta g(\bs{x};\varepsilon)$ denotes the variation in $g$ due to the topological change at $\bs{x}$. The topological change is characterised by the insertion of a rigid ball of radius $\varepsilon$ centred at $\bs{x}$. It is reasonable to think, according to the definition of the topological derivative \eqref{def:topological_derivative}, that the RHS of \equref{numerical_difference} approximates the topological derivative of $g$ when $\varepsilon$ is sufficiently small. Let us consider the scatterers and observation point as \figref{section7subsec1_configuration}. The incident wave is the plane wave along the $x_2$-axis of which the angular frequency is in $[2.5,3.5]$. The wave velocity is 1, and the amplitude is 1. The topological derivative is evaluated at the 31 points shown in \figref{section7subsec1_configuration}. Figure \ref{Df_pt1_466} shows the curve of \equref{Dfpade} and numerical difference in $f$ vs $\omega$ at $\bs{x}=(-6.0,-6.0)$ and $(-6.0,0.0)$. The $[4,4]$ Pad\'{e} approximation at $\omega_0=3.0$ is used in ``Indirect approx.'' and a small circular scatterer with the radius of $\varepsilon=0.01$ discretised with 100 boundary elements is used in ``Numerical difference''. In \figref{Df_pt1_466}, a pair of dashed lines show the range where the absolute error of ``Indirect approx.'' against ``Numerical difference'' is less than 0.01. We can see that ``Indirect approx.'' agrees with ``Numerical difference'' in some range of the angular frequency around $\omega_0$. The width of the range is, however, narrower than the range where $f(\omega)$ is accurately approximated (See \figref{f_plot}), which implies the necessity to monitor the accuracy of $\mathcal{D}_\mathrm{T}J$ rather than $J$ when we determine the number of subintervals (See Section \ref{range_division}).
The ranges where the error is less than 0.01 are plotted for all the evaluation points in \figref{range_comparison}. One finds that the range depends on the location where the topological derivative is evaluated. We, therefore, consider all the points where the computation of the topological derivative is performed when determining the number of subintervals (See Section \ref{range_division}). In the case of \figref{range_comparison}, the number of subintervals gets three in the case that the target band is $\omega\in[2.5, 3.5]$ based on Eqs. (\ref{omega_R_L_def1})--(\ref{omega_R_L_def2}) and Eqs. (\ref{define_N_div1})--(\ref{define_N_div2}). The topological derivative of the objective function $J$ evaluated by \equref{DJ_indirect} with the three subintervals is shown in \figref{DJ_plot} (``Indirect approx.''). For comparison, the corresponding numerical difference in $J$ by \equref{numerical_difference} is also plotted in \figref{DJ_plot} (``Numerical difference''). In ``Numerical difference'', the angular frequency integral in \equref{objective_function} is evaluated by the trapezoidal rule with the 100 integral points and the radius of the small circular scatterer is $\varepsilon=0.01$ in \equref{numerical_difference}. We can see that ``Indirect approx.'' agrees with ``Numerical difference'' in \figref{DJ_plot}.

\subsection{Topology optimisations}
In this subsection, we show the examples of the topology optimisation with our method.
\subsubsection{Acoustic lens}
\begin{figure}
  \centering
  \scalebox{0.9}{\includegraphics{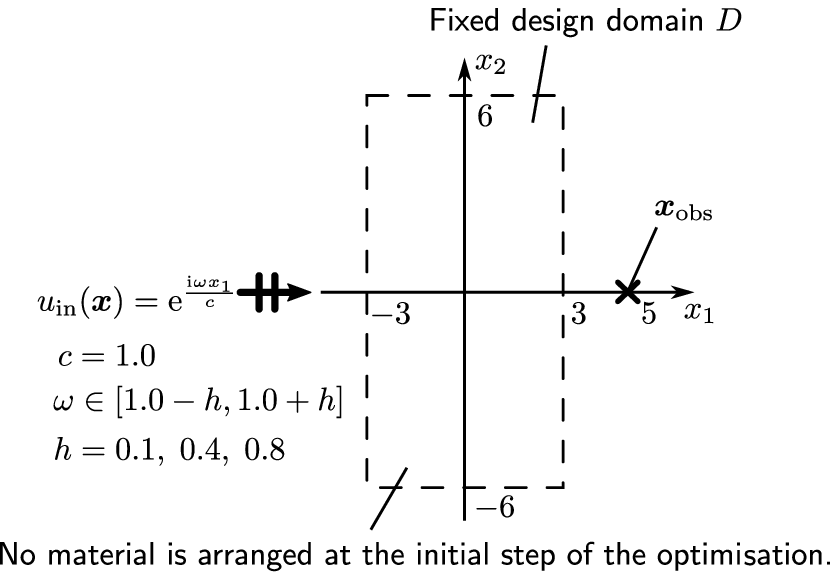}}
  \caption{Settings for the topology optimisation of acoustic lens.}
  \label{lens_settings}
\end{figure}
\begin{figure}
  \centering
  \scalebox{0.4}{\includegraphics{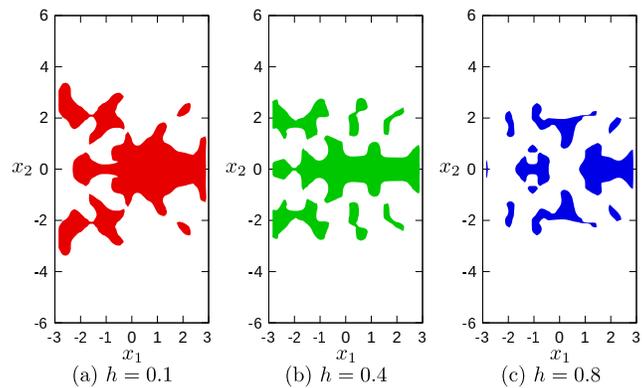}}
  \caption{The optimised acoustic lenses. The coloured area shows the rigid bodies.}
  \label{lens_mesh}
\end{figure}
\begin{figure}
  \centering
  \scalebox{0.4}{\includegraphics{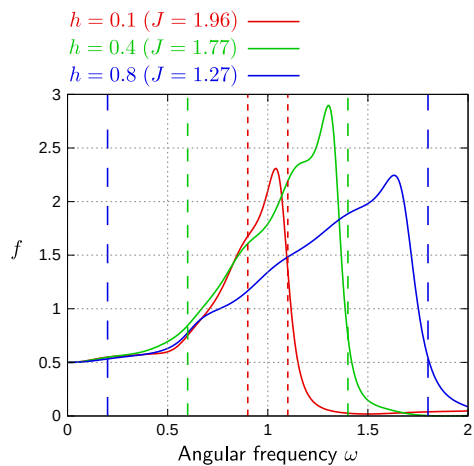}}
  \caption{The angular frequency responses of the optimised acoustic lenses with different target bands. A pair of the same coloured vertical dashed lines show the corresponding target band.}
  \label{lens_f_plots}
\end{figure}
First, we optimise acoustic lenses under the three conditions with different target bands. Figure \ref{lens_settings} shows the configuration of the optimisation. Three acoustic lenses are optimised to maximise the objective function \equref{objective_function} defined over the different target bands of which the midpoints are shared. We use the $[3,3]$ Pad\'{e} approximation $u_{[3,3]}$ for the frequency response estimation. The tolerance $\delta_{\mathcal{D}_\mathrm{T}f}$ in \equref{omega_R_L_def_3} is set as $0.01$, which is roughly $1\sim 10$ percent of the magnitude of $\mathcal{D}_\mathrm{T}f$.
Figures \ref{lens_mesh} and \ref{lens_f_plots} show the shapes of the three optimised acoustic lenses and their angular frequency responses $f(\omega)$. In \figref{lens_f_plots}, the acoustic lens optimised for a target band exhibits higher sound pressure in the band on average than the others, which implies the present topology optimisation successfully optimised the lenses. The acoustic lens optimised for $h=0.1$, for example, can focus more acoustic waves to the intended focal point $\bs{x}_\mathrm{obs}$ in $\omega\in[0.9,1.1]$. The angular frequency responses have their peaks in the higher frequency side in the bands and descend in the lower, which implies that the strategy to give priority to the performance in the higher frequency and suppress the deterioration of it in the lower is effective to improve the frequency average of the sound pressure.

\subsubsection{Sound shield}
\begin{figure}
  \centering
  \scalebox{1}{\includegraphics{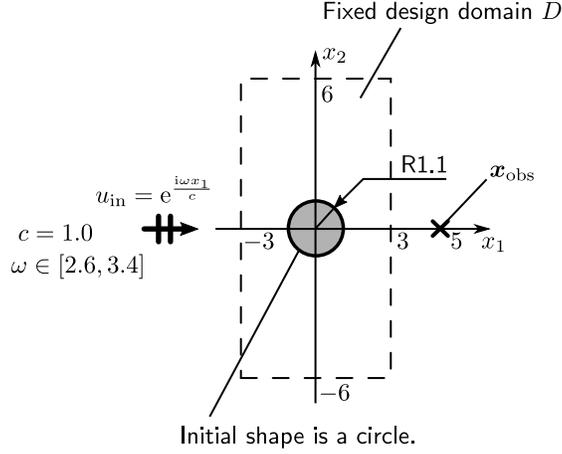}}
  \caption{Settings in the topology optimisation of sound shield.}
  \label{shield_settings}
\end{figure}
\begin{figure}
  \centering
  \scalebox{0.5}{\includegraphics{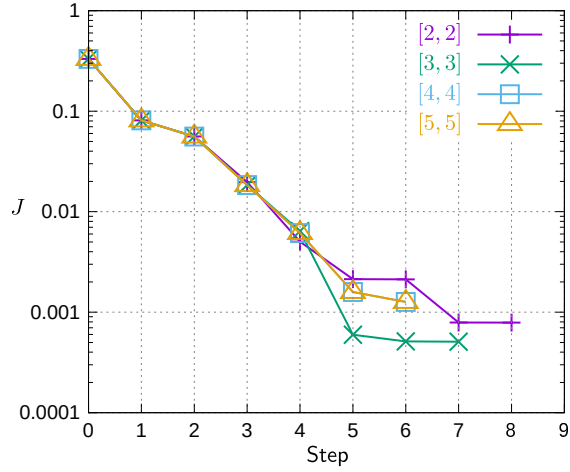}}\\
  (a) $\delta_{\mathcal{D}_\mathrm{T}f}=10^{-2}$\\
  \vspace{6mm}
  \scalebox{0.5}{\includegraphics{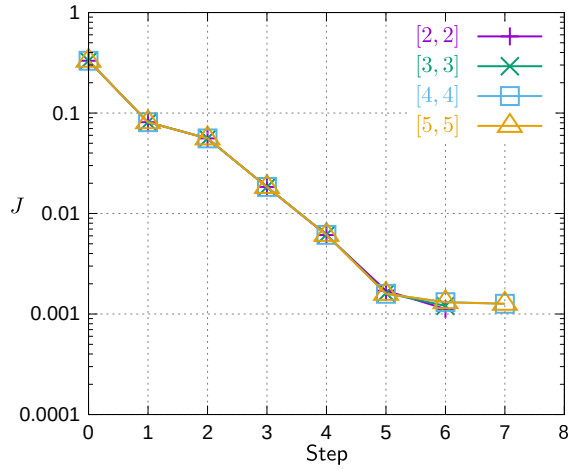}}\\
  (b) $\delta_{\mathcal{D}_\mathrm{T}f}=10^{-4}$
  \caption{The histories of the objective function in the optimisations when $[m,n]$ Pad\'{e} approximation is used.}
  \label{shield_J_histories}
\end{figure}
\begin{figure}
  \centering
  \scalebox{0.49}{\includegraphics{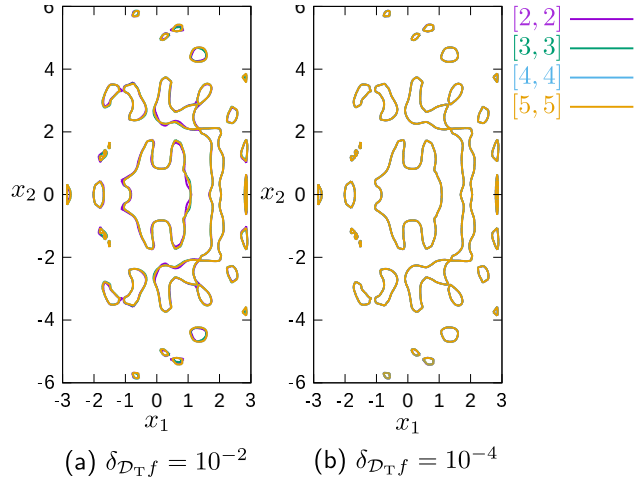}}\\
  \caption{The optimised structures of the sound shields. (a) The structures have minor differences depending on the degrees of the Pad\'{e} approximation if the tolerance for the frequency-band subdivision is large, (b) while the results with different degrees are consistent with each other if the tolerance is sufficiently small.}
  \label{shield_mesh}
\end{figure}
\begin{figure}
  \centering
  \scalebox{0.4}{\includegraphics{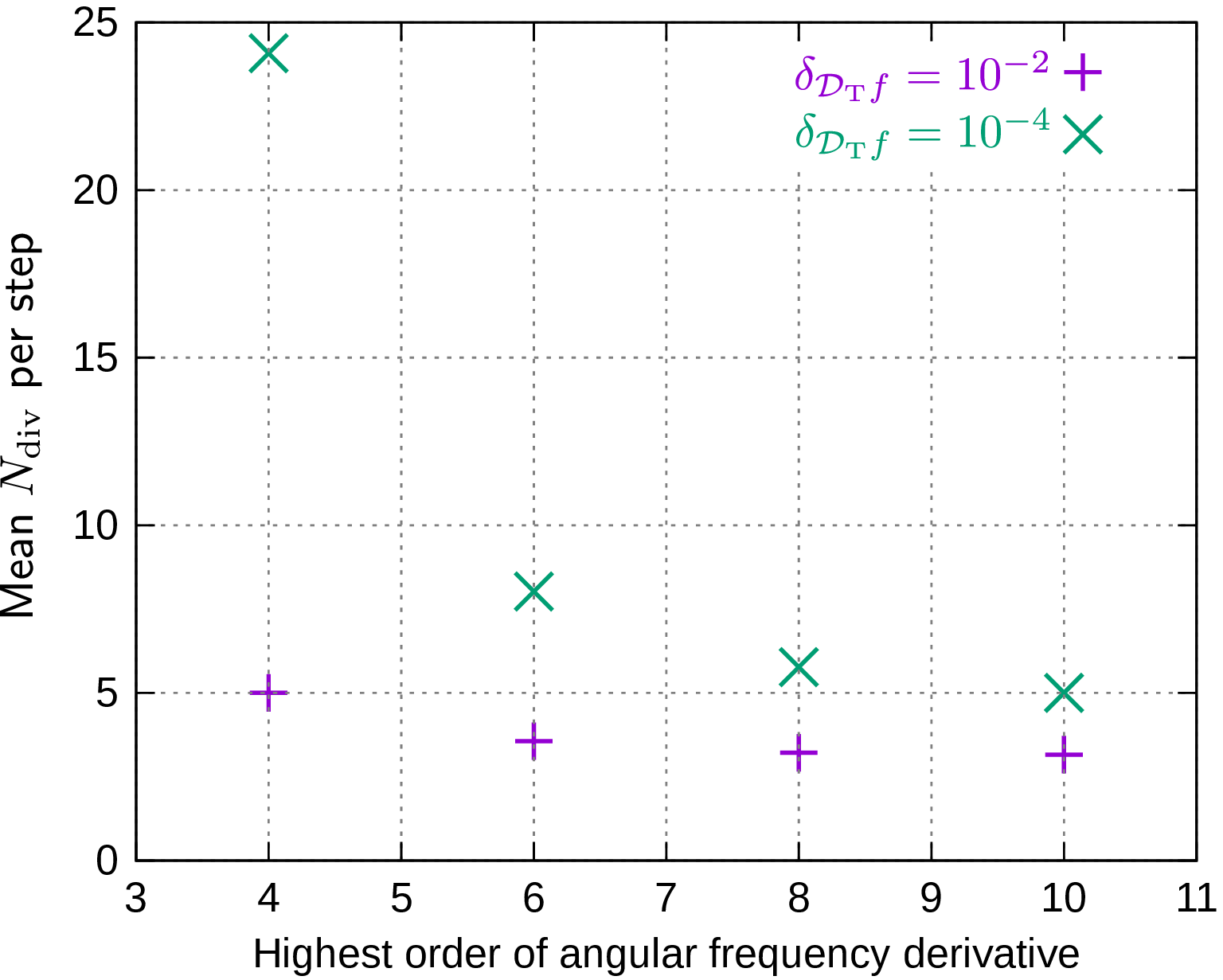}}
  \caption{The mean number of the subintervals per step.}
  \label{mean_Ndiv_per_step}
  \vspace{6mm}
  \scalebox{0.4}{\includegraphics{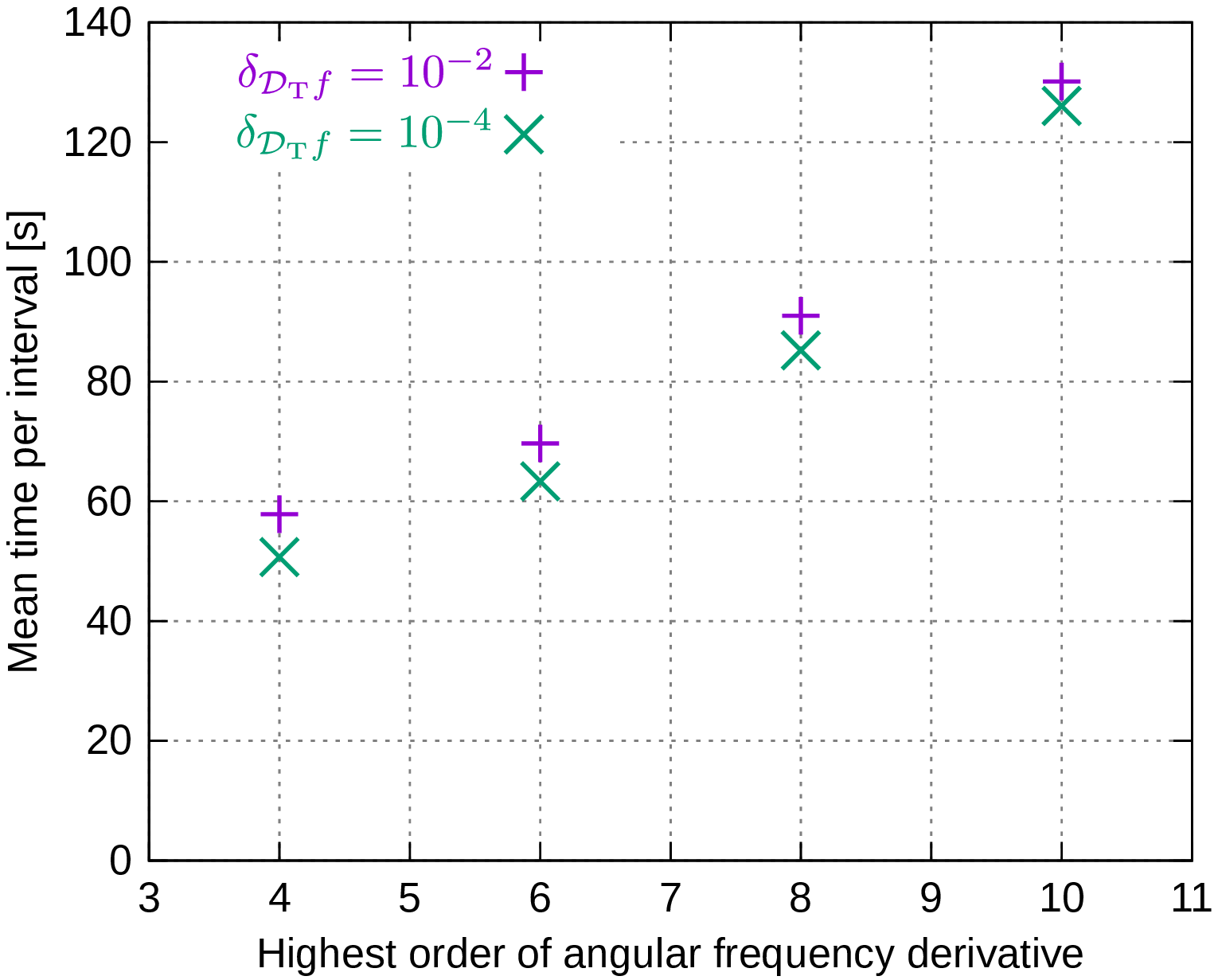}}
  \caption{The mean computational time per interval. The used CPU is Intel Xeon E5-2690 v3.}
  \label{mean_time_per_interval}
  \vspace{6mm}
  \scalebox{0.4}{\includegraphics{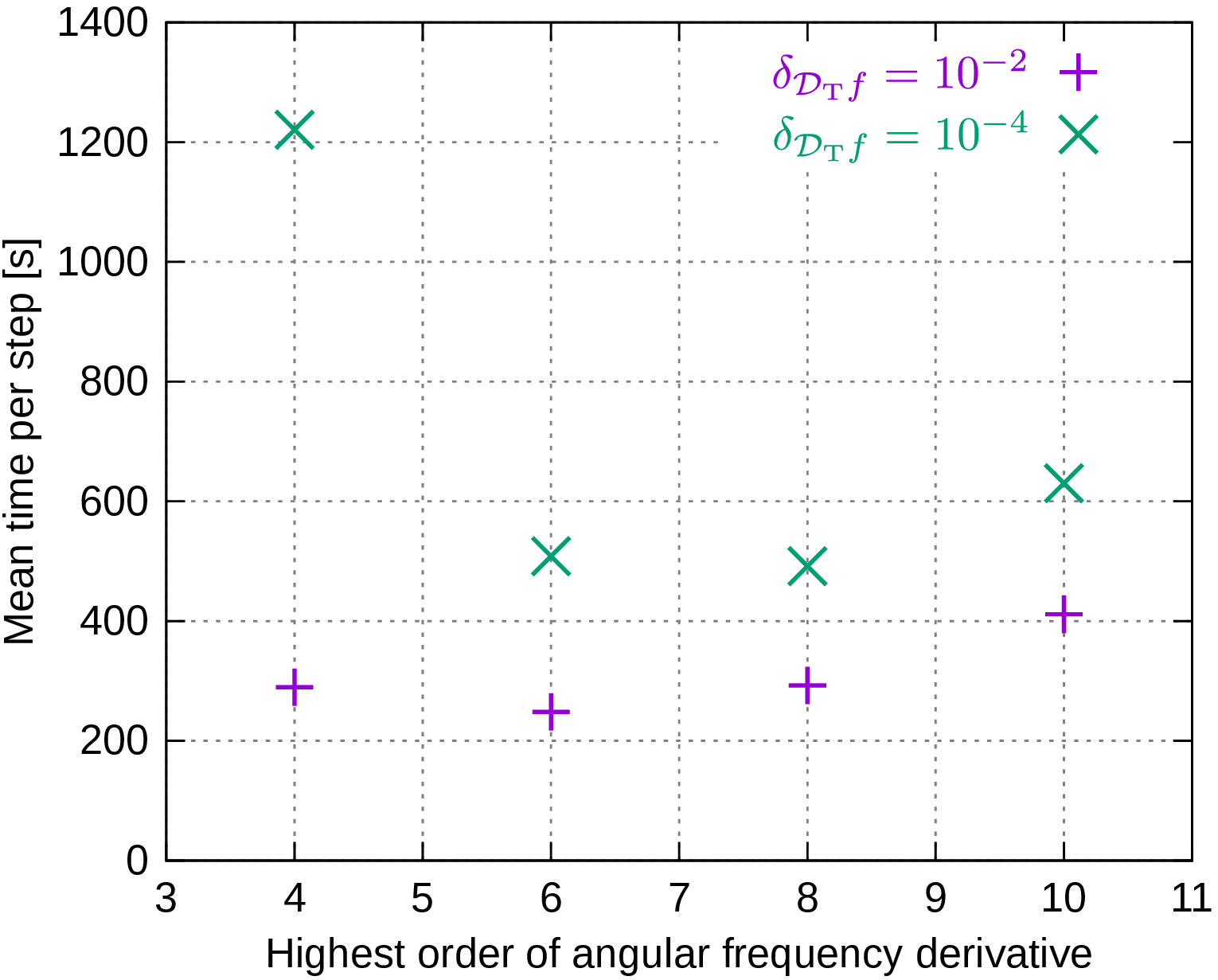}}
  \caption{The mean computational time per step. The used CPU is Intel Xeon E5-2690 v3.}
  \label{mean_time_per_step}
\end{figure}
Next, we optimise sound shields. Figure \ref{shield_settings} shows the configuration of the optimisation. To investigate the effects of the tolerance $\delta_{\mathcal{D}_\mathrm{T}f}$ in \equref{omega_R_L_def_3}, we optimise sound shields under the following two settings: $\delta_{\mathcal{D}_\mathrm{T}f}=10^{-2}$ and $10^{-4}$. In addition to $\delta_{\mathcal{D}_\mathrm{T}f}$, we also change the degrees of the Pad\'{e} approximation as $[2,2]$, $[3,3]$, $[4,4]$, and $[5,5]$.
Figures \ref{shield_J_histories} and \ref{shield_mesh} show the histories of the objective function in optimisation and the optimised sound shield structures. In \figref{shield_J_histories} (a), we can see that the optimisations follow a different path than the others when $[2,2]$ and $[3,3]$ Pad\'{e} approximation are used. This is because the error in the approximation accumulates as the optimisation proceeds. On the other hand, in \figref{shield_J_histories} (b), the optimisation histories are the same regardless of the degrees of the Pad\'{e} approximation.
We can infer that the smaller $\delta_{\mathcal{D}_\mathrm{T}f}$ leads to the more rigorous estimations of the objective function and topological derivative.

Figures \ref{mean_Ndiv_per_step}, \ref{mean_time_per_interval}, and \ref{mean_time_per_step} show the mean number of subintervals, the mean computational time per interval, and the mean computation time per step throughout an optimisation, respectively. The fewer subintervals are required when we use the angular frequency derivatives up to the higher-order as we can see in \figref{mean_Ndiv_per_step}. The fewer subintervals are desirable in terms of the computational cost because the fewer boundary value problems are solved. On the other hand, more computational cost is required to evaluate the higher-order derivatives as we can see in \figref{mean_time_per_interval}. Due to this trade-off, the computational cost gets the smallest with the angular frequency derivatives up to a specific order. In this sound shield optimisation, it is 6th ($[3,3]$ Pad\'{e} approx.) or 8th order ($[4,4]$ Pad\'{e} approx.) as we can see in \figref{mean_time_per_step}.

\section{Conclusion}\label{conclusion_section}
In this paper, we proposed a method to evaluate the frequency average of squared effective sound pressures and its topological derivative that uses the Pad\'{e} approximation and analytical integral. We also showed some numerical examples and confirmed the validity of our method. In Section \ref{objective_func_section}, we proposed two methods to approximate the angular frequency response; the direct approximation method and the indirect one. The response of the squared effective sound pressure is directly approximated in the former, while the response of the complex sound pressure is first approximated, and the squared effective sound pressure is then derived in the latter. We confirmed that the latter is better in the sense of the convergence speed through some numerical tests. In Section \ref{sensitivity_section}, we derived the topological derivative of the objective function based on the indirect approximation. First, the topological derivative of the complex sound pressure at the observation points is evaluated by the adjoint variable method. Then, the topological derivative of the objective function is derived by the chain rule. In Section \ref{range_division}, we proposed a simple method to determine the number of the target band subintervals. The numerical test in Section \ref{example_section} suggests the necessity to monitor the topological derivative accuracy when we determine the number of the subintervals because the frequency range where the approximation of the topological derivative holds is narrower than that of the objective function itself. We, therefore, proposed a method that estimates the frequency range where the approximation of the topological derivative is valid by using the lower-order approximation. The division of the target bandwidth by the estimated range width gives the number of subintervals. In Section \ref{example_section}, we showed some numerical examples. The optimisation examples of the acoustic lens show that by giving the different target bands, we can obtain the acoustic lens optimised to each band as expected. Through the optimisation examples of sound shields, we checked the effect of the tolerance that controls the approximation accuracy. It is confirmed that we can conduct a more rigorous optimisation by imposing a more strict tolerance. In addition, we found that we can minimise the computational cost by choosing the proper degrees of the Pad\'{e} approximation.

In this paper, we showed only the simple optimisation examples; acoustic lens and sound shield. On the other hand, it is known that sound pressure is more sensitive to the frequency in periodic structures and wave-guide structures. The superiority of the proposed method; the integration strictness may be more important in the optimisation of such structures because the function with sharp peaks and dips is generally hard to integrate numerically. In addition, our objective function does not flatten the optimised frequency response, unlike the robust topology optimisation. If the frequency response should be flat in applications, we need to use different objective functions. One option is to combine the robust topology optimisation and the Pad\'{e} approximation. Another simple option is to make the frequency response function $f(\omega)$ flat by minimising $\left|\frac{\mathrm{d}f}{\mathrm{d}\omega}\right|$ in the target band. One of the other possible improvements is the target band subdivision. In this paper, we divide the target band into equally spaced subintervals for simplicity. By changing the interval length adaptively in the target band, we may improve the approximation accuracy with a lower computational cost.

% \section*{References}
% \bibliography{ref}

\appendix

\section{Level-set method with B-spline surface}
  \label{level_set_section}
In this section, a topology optimisation method with the level-set of a B-spline surface \cite{isakari2017Bspline} is reviewed, which is used to solve the structural optimisation problem. The use of this method is not essential, and other solvers such as \cite{yamada2010topology} are also available.

In the level-set method, a scalar function $\phi(\bs{x})$ called the level-set function is introduced to represent the material distribution in the design domain $D$, which is defined as
 \begin{align}
   &
   \begin{cases}
     \phi(\bs{x})>0 & \bs{x}\in\Omega,\\
     \phi(\bs{x})=0 & \bs{x}\in\Gamma,\\
     \phi(\bs{x})<0 & \bs{x}\in D\backslash\overline{\Omega},
   \end{cases}\\
   \label{level_set_norm}
   &\quad||\phi(\bs{x})||=1,
 \end{align}
where $||\cdot||$ denotes a norm defined on $D$. Amstutz and Andr\"a~\cite{amstutz2006new} proposed to update the level-set function by
 \begin{equation}
   \frac{\partial\phi(\bs{x},t)}{\partial t}=[\mathcal{D}_\mathrm{T}J](\bs{x},t)-(\mathcal{D}_\mathrm{T}J,\phi)\phi(\bs{x},t),
 \end{equation}
where $t$ and $(\cdot,\cdot)$ denote the fictitious time corresponding to the optimisation step and the inner product that defines the norm in \equref{level_set_norm}, respectively. In the method proposed in \cite{isakari2017Bspline}, the level-set function $\phi(\bs{x})$ is discretised in space by the B-spline surface, which gives a way to control the geometric complexity of the optimised structure by adjusting the degree of freedom in the B-spline surface. In addition to $\phi(\bs{x})$, the topological derivative $[\mathcal{D}_\mathrm{T}J](\bs{x})$ is also discretised by the B-spline surface which shares the same basis functions and knots with $\phi(\bs{x})$.
In the original paper \cite{isakari2017Bspline}, the control variables of the B-spline surface are determined such that the surface interpolates $\mathcal{D}_\mathrm{T}J$ at the Greville abscissae in the design domain $D$. This method may, however, give the optimised structure with wavy boundaries due to the Gibbs-like phenomenon. To avoid this problem, we use an improved version of the method, in which the control variables are determined such that the $H_1$ norm discrepancy between the B-spline surface and $\mathcal{D}_\mathrm{T}J(\bs{x})$ is minimised. The norm in \equref{level_set_norm} and corresponding inner product is also defined in the sense of $H_1$ accordingly.

\section{The Pad\'{e} approximation}
  \label{pade_section}
  In this section, general features of the Pad\'{e} approximation are briefly described. The Pad\'{e} approximation is a method to approximate a function $f(x)$ with the rational polynomial as:
  \begin{equation}
    \label{Rx}
    f(x)\simeq f_{[M,N]}(x):=\frac{p(x)}{q(x)}:=\frac{\displaystyle{\sum_{i=0}^{M}}p_i(x-x_0)^i}{\displaystyle{\sum_{i=0}^{N}}q_i(x-x_0)^i},
  \end{equation}
  which is called the $[M,N]$ Pad\'{e} approximation. The coefficients $p_i$ and $q_i$ are unique up to constant and are generally normalised as $q_0=1$. The rest of the coefficients are determined such that the Taylor series of $f_{[M,N]}(x)$ is identical to that of $f(x)$ up to the $(M+N)$-th order. Let us denote the Taylor series of $f(x)$ as
  \begin{equation}
    \label{fx}
    f(x)=\sum_{i=0}^{\infty}a_i(x-x_0)^i,
  \end{equation}
  with $a_i:=f^{(i)}(x_0)/{i!}$. Multiplying $q(x)$ to Eqs. \eqref{Rx} and \eqref{fx}, we have
  \begin{align}
    \label{Rxqx}
    f_{[M,N]}(x)q(x)&=\sum_{i=0}^{M}p_i(x-x_0)^i,\\
    \label{fxqx}
    f(x)q(x)&=\left[\sum_{i=0}^{\infty}a_i(x-x_0)^i\right]\left[\sum_{i=0}^{N}q_i(x-x_0)^i\right],
  \end{align}
  respectively. The coefficients $p_i$ and $q_i$ are determined such that Eqs. \eqref{Rxqx} and \eqref{fxqx} are equal to each other up to the $(M+N)$-th order, which is equivalent to the previous definition. Expanding \equref{fxqx} and comparing it with \equref{Rxqx}, we have the following $M+N+1$ equations:
  % \begin{equation*}
  %   \begin{aligned}
  %     p_0&=a_0\\
  %     p_1&=a_1+a_0q_1\\
  %     p_2&=a_2+a_1q_1+a_0q_2\\
  %     &\hspace{1.5mm}\vdots\\
  %     p_i&=a_i+\sum_{j=1}^ia_{i-j}q_j\quad(1\leq i\leq N)\\
  %     &\hspace{1.5mm}\vdots\\
  %     p_n&=a_n+\sum_{j=1}^N a_{N-j}q_j\\
  %     p_{N+1}&=a_{N+1}+\sum_{j=1}^N a_{N+1-j}q_j\\
  %     &\hspace{1.5mm}\vdots\\
  %     p_i&=a_i+\sum_{j=1}^N a_{i-j}q_j\quad(N<i\leq M)\\
  %     &\hspace{1.5mm}\vdots\\
  %     p_m&=a_m+\sum_{j=1}^N a_{M-j}q_j\\
  %     0&=a_{M+1}+\sum_{j=1}^N a_{M+1-j}q_j\\
  %     &\hspace{1.5mm}\vdots\\
  %     0&=a_i+\sum_{j=1}^N a_{i-j}q_j\quad(M<i\leq M+N)\\
  %     &\hspace{1.5mm}\vdots\\
  %     0&=a_{M+N}+\sum_{j=1}^N a_{M+N-j}q_j
  %   \end{aligned}
  % \end{equation*}
  % which can be compressed as
  \begin{align}
    p_0&=a_0, \\
    \label{Pade_tobe_solved}
    p_i&=a_i+\sum_{j=1}^ia_{i-j}q_j\quad(1\leq i\leq M+N),
  \end{align}
with $p_i:=0\;(i>M)$ and $q_i:=0\;(i>N)$ defined.% The same result is obtained when $M<N$, too.

  In this study, the linear equations \eqref{Pade_tobe_solved} are solved by the GMRES (Generalised Minimal RESidual) \cite{saad1986gmres} without the restart because $M+N$ is at most 10 or 20, which guarantees the numerical solution converges to the exact one up to the rounding error after the $M+N$ iterations. The iterative solver is preferable to the direct solver because \equref{Pade_tobe_solved} may get singular when $f(x)$ itself is a rational function and the degrees of the Pad\'{e} approximation surpass those of $f(x)$ . When $f(x)$ is, for example, given by
\begin{equation}
  f(x)=\frac{1}{(x-1)(x-3)}
\end{equation}
and is approximated by the $[5,5]$ Pad\'{e} approximation at $x=2$, we have
\begin{align}
  f_{[M,N]}(x)&=\frac{-1-(x-2)^2/2}{1-(x-2)^2/2-(x-2)^4/2}\\
  \label{canceled_pair}
  &=\frac{(x-2+\mathrm{i}\sqrt{2})(x-2-\mathrm{i}\sqrt{2})}{(x-1)(x-3)(x-2+\mathrm{i}\sqrt{2})(x-2-\mathrm{i}\sqrt{2})},
\end{align}
in which the factors $x-2+\mathrm{i}\sqrt{2}$ and $x-2-\mathrm{i}\sqrt{2}$ are canceled. This means that the corresponding terms are not necessarily the ones in \equref{canceled_pair}, i.e. the coefficients of $f_{[M,N]}$ cannot be determined uniquely. In such a case, the coefficient matrix in the algebraic equation \eqref{Pade_tobe_solved} is singular. Although we have not encountered any singular matrix throughout our experiments, we recommend using the GMRES to solve \equref{Pade_tobe_solved} just in case.

\end{document}